\documentclass[11pt,a4paper,oneside]{amsart}
\usepackage[margin=25mm,headheight=14pt,headsep=7mm]{geometry}
\usepackage[T1]{fontenc}
\usepackage{lmodern}
\usepackage{amssymb,mathrsfs,comment,braket}
\usepackage{needspace}
\usepackage[colorlinks=true,linkcolor=black,citecolor=black,urlcolor=blue,
pdfauthor={Dongwen Liu and Dipendra Prasad},
pdftitle={Restricting representations of finite and p-adic reductive groups to unipotent subgroups}]{hyperref}
\newcommand{\gr}{\operatorname{gr}}

\newcommand{\St}{\operatorname{St}}
\newcommand{\Ind}{\operatorname{Ind}}
\newcommand{\cind}{\operatorname{c-Ind}}
\newcommand{\GL}{\mathrm{GL}}

\newcommand{\Cc}{C_c^\infty}
\newcommand{\C}{\mathbb C}
\newcommand{\one}{\mathbf1}
\newtheorem{theorem}{Theorem}[section]
\newtheorem{conjecture}[theorem]{Conjecture}
\newtheorem{lemma}[theorem]{Lemma}
\newtheorem{proposition}[theorem]{Proposition}
\newtheorem{corollary}[theorem]{Corollary}
\newtheorem{example}[theorem]{Example}

\theoremstyle{remark}
\newtheorem{remark}[theorem]{Remark}
\numberwithin{equation}{section}
\title[Restrictions to unipotent subgroups]{Restricting representations of finite and p-adic reductive groups to unipotent subgroups}

\author[D. Liu]{Dongwen Liu}
\address{School of Mathematical Sciences, Zhejiang University, 866 Yuhangtang Road, Hangzhou 310058, P.R. China}
\email{maliu@zju.edu.cn}

\author[D. Prasad]{Dipendra Prasad}
\address{Department of Mathematics, Indian Institute of Technology Bombay, Powai, Mumbai 400076, India}
\email{prasad.dipendra@gmail.com}

\date{29 November 2026, 9:54 p.m. IST}

\begin{document}

\begin{abstract}
We formulate a conjecture about nonvanishing of $\pi_{N,\psi}$, where $\psi:N\to\C^\times$ is any character of any unipotent subgroup $N$ of a quasi-split connected reductive group $G$ over a finite or nonarchimedean local field, and $\pi$ is any irreducible generic representation of $G$. In characteristic zero or sufficiently large positive characteristic, we prove the conjecture for $\GL_n$ and for the Steinberg representation of any connected reductive group. We also consider the restriction of any irreducible representation of $G$ to a maximal unipotent subgroup $U$, or to a maximal $F$-split connected solvable subgroup $B_0=AU$, where $A$ is a maximal $F$-split torus of $G$ normalizing $U$. We formulate two conjectures about these restrictions in appropriate Grothendieck groups, for stable representations over local fields and uniform representations over finite fields.
\end{abstract}

\maketitle

\setcounter{tocdepth}{1}

\tableofcontents

\section{Introduction}\label{sec:introduction}

Let $p$ be a prime number, and let $F$ be a nonarchimedean local field of residual characteristic $p$ or a finite field of characteristic $p$, with a fixed nontrivial additive character $\psi_F:F\to\C^\times$. In this paper, finite groups are endowed with the discrete topology, and representations of groups over nonarchimedean local fields or finite fields are over the complex numbers and are assumed to be smooth.

Let $\mathbf G$ be a connected quasi-split reductive group over $F$, and let $G:=\mathbf G(F)$. For a representation $\pi$ of $G$, a closed subgroup $N\subset G$ and a character $\psi:N\to\C^\times$, put
\[
\pi_{N,\psi}:=\pi/\operatorname{Span}\{\pi(u)v-\psi(u)v:u\in N,\ v\in\pi\}.
\]
Then
\[
\pi_{N,\psi}\ne\{0\}\quad\text{if and only if}\quad\operatorname{Hom}_N(\pi,\psi)\ne\{0\},
\]
where $\operatorname{Hom}_N(\pi,\psi)$ denotes the space of $(N,\psi)$-equivariant continuous linear functionals on $\pi$. When $\psi$ is trivial, we simply write $\pi_N=\pi_{N,\psi}$.

For $a\in\mathbf G_{\mathrm{ad}}(F)$, where $\mathbf G_{\mathrm{ad}}$ is the adjoint group of $\mathbf G$, we write $c_a$ for the corresponding automorphism of $G$ induced by conjugation, and define a representation $\pi^a$ of $G$ by
\[
\pi^a(g):=\pi(c_a^{-1}(g)),\qquad g\in G.
\]
For generic representations, we formulate the following conjecture on the nonvanishing of $\pi_{N,\psi}$.

\begin{conjecture}\label{conj:degenerate}
There is a constant $n_{\mathbf G}>0$, depending only on the absolute root datum of $\mathbf G$, with the following property. Suppose that $\operatorname{char}F=0$ or $\operatorname{char}F=p\ge n_{\mathbf G}$. Let $\pi$ be an irreducible representation of $G$ which is generic for some Whittaker datum, let $\mathbf N$ be a smooth connected unipotent $F$-subgroup of $\mathbf G$, and let $\psi$ be a character of $N=\mathbf N(F)$. Then there exists $a\in\mathbf G_{\mathrm{ad}}(F)$ such that
\[
(\pi^a)_{N,\psi}\ne\{0\},
\]
unless some proper $F$-parabolic subgroup of $\mathbf G$ has unipotent radical $\mathbf N'$ for which
\[
N':=\mathbf N'(F)\subseteq N,\qquad \psi|_{N'}=1,\qquad \pi_{N'}=\{0\}.
\]
\end{conjecture}

\begin{remark}
If a representation $\pi$ of $G$ contains $\Cc(U)$ upon restriction to a maximal unipotent subgroup $U$, as is the case for a principal series representation arising from a minimal parabolic subgroup of $\mathbf G$ and, as we will see in Section~\ref{sec:steinberg}, also for the Steinberg representation, then clearly Conjecture~\ref{conj:degenerate} holds.
\end{remark}

\begin{remark}\label{rem:adjoint-group}
If $\mathbf G$ and $\mathbf G'$ are two reductive groups over $F$ with the same adjoint group, then Conjecture~\ref{conj:degenerate} holds for $\mathbf G$ if and only if it holds for $\mathbf G'$.
\end{remark}

We prove Conjecture~\ref{conj:degenerate} for $\mathbf G=\GL_n$, with $n\ge1$, in the following more precise form. Its proof is given in Theorem~\ref{thm:precise}. Over a finite field of characteristic $p\ge n$, $N$ may be any $p$-subgroup.

\begin{theorem}\label{thm:main}
Let $G=\GL_n(F)$, where $\operatorname{char}F=0$ or
$\operatorname{char}F=p\ge n$, and let $\pi$ be an irreducible generic
representation of $G$. If $F$ is local, let $N=\mathbf N(F)$ for a
smooth connected unipotent subgroup $\mathbf N$ of $\GL_n$; if $F$ is
finite, let $N$ be any $p$-subgroup of $G$. Let $\psi$ be a character of $N$.

Conjugating if necessary, assume that $N$ is contained in the upper
triangular maximal unipotent subgroup. Let $N_i$ denote the unipotent radical of the standard maximal parabolic
of type $(i,n-i)$, for $1\leq i<n$. Let $P$ be the standard parabolic
whose unipotent radical contains exactly those $N_i\subseteq N$ for which
$\psi|_{N_i}=1$. Then $\pi_{N,\psi}=0$ if and only if the Jacquet
module of $\pi$ along $P$ is zero. If there are no such $N_i$, then
$\pi_{N,\psi}\ne 0$.
\end{theorem}

\begin{example}[Containment on rational points]\label{ex:rational-containment}
In the upper unitriangular subgroup of $\GL_3$ over $\mathbb F_q$, let $\mathbf N$ be defined by $u_{23}=u_{12}^q-u_{12}$. Since $x^q=x$ on $\mathbb F_q$, $\mathbf N(\mathbb F_q)=\mathbf V(\mathbb F_q)$, where $\mathbf V$ is the unipotent radical of the standard $(1,2)$-parabolic. Thus $\pi_{\mathbf N(\mathbb F_q)}=0$ for every cuspidal $\pi$, which is generic, and $\pi^a\simeq\pi$ for all $a\in\mathrm{PGL}_3(\mathbb F_q)$. Yet $\mathbf N$ contains no proper parabolic radical. Indeed, $\mathbf N$ is an affine plane with coordinates $u_{12},u_{13}$, and the $(1,3)$-entry of a commutator is $xy^q-yx^q$ for $(1,2)$-entries $x,y$. Hence $\mathbf N$ is nonabelian over $\overline{\mathbb F}_q$, while the proper parabolic radicals in $\GL_3$ are either two-dimensional and abelian or three-dimensional. This explains the containment requirement on $F$-points (rather than over $\overline F$) in Conjecture~\ref{conj:degenerate}.
\end{example}

\begin{example}[The characteristic bound]\label{ex:characteristic-bound}
Some restriction on the characteristic is necessary in Theorem~\ref{thm:main}.
Let $G=\GL_3(\mathbb F_2)$ and
\[
 N=\left\{\begin{pmatrix}1&a&b\\0&1&a\\0&0&1\end{pmatrix}:
 a,b\in\mathbb F_2\right\}=\langle r\rangle\simeq\mathbb Z/4\mathbb Z,
 \qquad r=\begin{pmatrix}1&1&0\\0&1&1\\0&0&1\end{pmatrix}.
\]
There are two three-dimensional irreducible representations $\pi$ of $G$.
Both are cuspidal, hence generic, and their restrictions to $N$ have character values $(3,1,-1,1)$
on $(1,r,r^2,r^3)$, so $\pi_{N,\psi}=0$ for $\psi(r)=-1$.
In this case, $\ker\psi$ has order two, whereas every proper parabolic
subgroup has a unipotent radical of order at least four. Therefore,
the unipotent radical of no proper parabolic subgroup is contained
in $\ker\psi$. Moreover, $N$ is the group of rational points of
the smooth connected $\mathbb F_2$-split unipotent subgroup defined by
the same matrix formula over $\mathbb F_2$.
\end{example}

Let us briefly explain the proof of Theorem~\ref{thm:main}. If the
ordinary Jacquet module specified in the theorem is zero, so is
$\pi_{N,\psi}$, by taking coinvariants in stages. Otherwise we pass to
an irreducible generic quotient of that Jacquet module and work block
by block in the Levi subgroup. In a single block the idea is to use the
bottom piece of the Bernstein--Zelevinsky filtration and find a maximal
unipotent subgroup $U$ in the mirabolic subgroup, together with a
nondegenerate character $\theta:U\to\C^\times$ compatible with $(N,\psi)$
in the sense that
\[
\theta|_{U\cap N} = \psi|_{U\cap N}.
\] 
Since the bottom piece of the Bernstein--Zelevinsky filtration is compactly induced from the character $\theta$ of $U$ to the mirabolic subgroup, Mackey theory, combined with the exactness of the functor $\pi\mapsto\pi_{N,\psi}$, gives the \emph{nonvanishing} of $\pi_{N,\psi}$. A topological point is needed here: the image of $N$ in $U\backslash P_n$ must be closed (equivalently, $UN$ must be closed in $P_n$). This is not automatic for closed subgroups of a topological group; it follows from Lemma~\ref{lem:closed-orbit}. In positive characteristic, the assumption $p\geq n$ is needed to construct the compatible pair $(U,\theta)$: it allows us to use the exponential map and the fact that fewer than $p$ proper affine subspaces cannot cover an affine space over $\mathbb F_p$. See the proof of Proposition~\ref{prop:extension} for details.

For generalized Whittaker models, a Jordan basis for the nilpotent element gives a direct construction of compatible characters. Corollary~\ref{cor:generalized} thus shows that every irreducible generic representation of $\GL_n(F)$ has all these models over any finite or nonarchimedean local field $F$, including fields of characteristic $p<n$. The following example corresponds to the partition $[2^r]$. For convenience, when there is no confusion we also identify an algebraic group over $F$ with the group of its $F$-points, e.g.\ a unipotent subgroup of $\mathbf G$ refers to the group of $F$-points of a unipotent subgroup of $\mathbf G$.

\begin{example}
Let $F$ be any finite or nonarchimedean local field and put $G=\GL_{2r}(F)$. Let $N$ be the unipotent radical of the standard maximal parabolic subgroup of $\GL_{2r}(F)$ of type $(r,r)$, and define a character
\[
\psi:N\longrightarrow\C^\times,\qquad
\begin{pmatrix}1_r&X\\0&1_r\end{pmatrix}\longmapsto\psi_F(\operatorname{tr}X).
\]
If $\pi$ is an irreducible generic representation of $\GL_{2r}(F)$, then
\[
\pi_{N,\psi}\ne\{0\}.
\]
\end{example}

Extending our study of characters of unipotent subgroups of $G=\mathbf G(F)$, with $\mathbf G$ reductive, occurring in generic representations of $G$, we consider the restriction of \emph{any irreducible representation} of $G$ to a maximal unipotent subgroup $U$ or, more generally,
to a maximal $F$-split connected solvable subgroup $B_0$, i.e., $B_0=AU$, where $A$ is a maximal $F$-split torus of
$\mathbf G$ normalizing $U$. We will discuss this restriction problem for $B_0$
and $U$ and formulate a general restriction conjecture for each of them
(see Section~\ref{sec:question}).

Beyond the study of $\pi_{N,\psi}$ for
degenerate Whittaker models in the pioneering works of Kawanaka on
generalized Gelfand--Graev representations~\cite{K85}, Lusztig on unipotent
support~\cite{L92}, and M\oe glin--Waldspurger for $p$-adic
groups~\cite{MW87}, which have considerably influenced the subject,
the more general restriction problem that we consider, to our
knowledge, has not been considered in the literature.

This paper is organized as follows. Sections~\ref{sec:cyclic} and~\ref{sec:extension} contain some preliminary linear algebraic results that are used in Section~\ref{sec:main-proof} to prove Theorem \ref{thm:main}. In Section~\ref{sec:generalized} we prove the nonvanishing of all generalized Whittaker models for an irreducible generic representation of $\GL_n(F)$. Section~\ref{sec:steinberg} recalls the parabolic resolution and the open-cell realization of the Steinberg representation, due over nonarchimedean local fields to Borel--Serre \cite[\S3]{BS76}. We use the resulting identification $\St_G|_U\simeq\Cc(U)$ to prove Conjecture~\ref{conj:degenerate} for the Steinberg representation under the hypotheses of Corollary~\ref{cor:steinberg-conjecture}. Section~\ref{sec:grothendieck-B} introduces the Grothendieck groups of representations of $B$, $B_0$, and $U$ used below. Sections \ref{sec:principal} and \ref{sec:GL} treat principal series, and all irreducible smooth representations of $\GL_n(F)$, respectively. Section~\ref{sec:question} formulates two restriction conjectures which, roughly speaking, suggest that in appropriate Grothendieck groups, the restrictions to $B_0$ and to $U$ resemble those of principal series representations of $G$, except for the possible occurrence of signs in virtual linear combinations. In Section~\ref{sec:restriction-evidence}, we prove that these conjectures are ``hereditary'' under parabolic induction and reduce the finite-field problem to elliptic Green functions.

\medskip
\noindent\textbf{Acknowledgement.}
Conjecture~\ref{conj:degenerate} of this paper has been proposed by D.P.\ on several occasions over the past few years, including at the recent conference \emph{Automorphic representations, L-functions and Harmonic analysis} at BIMSA in June 2026. The case of general linear groups was proved when both authors were attending the conference \emph{Wavefront Sets of Representations and Endoscopic Classification} at BIRS, Hangzhou. The authors are grateful for the hospitality of the organizers and the working environments provided by these institutes. D.P.\ would like to thank Taiwang Deng and Santosh Pattanayak for some discussions.

D.L.\ was partially supported by National Natural Science Foundation of China No.~12526208.

\section{A lemma on cyclic vectors}\label{sec:cyclic}
In this section, $F$ denotes an arbitrary field. Fix a positive integer $r$, and let $F^r$ (resp.\ $F_r$) be the space of $r\times1$ (resp.\ $1\times r$) matrices over $F$. We identify $F_r$ with the linear dual of $F^r$ by
\[
\lambda(v):=\lambda\cdot v,\qquad\lambda\in F_r,\quad v\in F^r.
\]
Let $e_1,e_2,\ldots,e_r$ be the standard basis of $F^r$. Then we have the complete flag
\[
\mathcal F_\bullet:\{0\}\subset\mathcal F_1\subset\mathcal F_2\subset\cdots\subset\mathcal F_r,
\]
where
\begin{equation}\label{eq:standard-flag}
\mathcal F_k:=\operatorname{Span}_F\{e_1,e_2,\ldots,e_k\},\qquad k=1,2,\ldots,r.
\end{equation}
For $X\in\mathfrak{gl}_r(F)=\operatorname{End}_F(F^r)$, an element $\lambda\in F_r$ is called $X$-cyclic if
\[
\lambda,\lambda X,\ldots,\lambda X^{r-1}
\]
form a basis of $F_r$. Let $\mathfrak b_r(F)$ be the space of $r\times r$ upper triangular matrices over $F$. Then we have the following result.

\begin{lemma}\label{lem:cyclic}
Let $X\in\mathfrak{gl}_r(F)$ and $\lambda\in F_r$. There exists $Y\in\mathfrak b_r(F)$ such that $\lambda$ is $(X+Y)$-cyclic if and only if there is no $k\in\{1,2,\ldots,r\}$ such that
\begin{equation}\label{eq:cyclic-obstruction}
X(\mathcal F_k)\subset\mathcal F_k\quad\text{and}\quad \lambda|_{\mathcal F_k}=0.
\end{equation}
\end{lemma}
\begin{proof}
When $F$ is finite, the lemma easily follows from \cite[Theorem~1.1, Corollary~1.2]{GRS92}. In the sequel, assume that $F$ is infinite.

If \eqref{eq:cyclic-obstruction} holds for some $k\in\{1,2,\ldots,r\}$, then $\mathcal F_k$ is preserved by $X+Y$ for every $Y\in\mathfrak b_r(F)$. It follows that each $\lambda(X+Y)^i$ vanishes on $\mathcal F_k$, $i=0,1,\ldots,r-1$, hence $\lambda$ cannot be $(X+Y)$-cyclic.

Conversely, assume that there is no such $k$. We need to show that there exists $Y\in\mathfrak b_r(F)$ such that $\lambda$ is $(X+Y)$-cyclic, i.e.\ the polynomial function
\[
Y\longmapsto\det\begin{pmatrix}
\lambda\\ \lambda(X+Y)\\ \vdots\\ \lambda(X+Y)^{r-1}
\end{pmatrix}
\]
on $\mathfrak b_r(F)$ is not identically zero. Since $F$ is infinite, we may pass to any extension of $F$ and therefore assume that $F$ is algebraically closed in the rest of the proof.

For $k=1,2,\ldots,r$, put
\[
\mathcal E_k:=\{v\in\mathcal F_k:X(v)\in\mathcal F_k,\ \lambda(v)=0\}\subsetneq\mathcal F_k.
\]
Then
\begin{equation}\label{eq:Ek-dimension}
\dim\mathcal E_k<k.
\end{equation}
For every $A\in\mathfrak{gl}_r(F)$, the element $\lambda\in F_r$ is not $A$-cyclic if and only if $A$ has an eigenvector in $\ker(\lambda)$. Indeed, consider
\[
V_A:=\bigcap_{i=0}^{r-1}\ker(\lambda A^i).
\]
By Cayley--Hamilton, $V_A$ is $A$-invariant. If $\lambda$ is not $A$-cyclic, then $V_A\ne\{0\}$ hence $A$ has an eigenvector in $V_A$. Conversely, an eigenvector of $A$ in $\ker(\lambda)$ lies in $V_A$.

Suppose that $Y\in\mathfrak b_r(F)$ and $A=X+Y$ has an eigenvector
$v=\left[\begin{smallmatrix}v_1\\v_2\\\vdots\\v_r\end{smallmatrix}\right]\in\ker(\lambda)$, with eigenvalue $a\in F$. We normalize $v$ such that the last nonzero coordinate of $v$ is $v_k=1$, where $k\in\{1,2,\ldots,r\}$. Then $v\in\mathcal F_k$, $Y(v)\in\mathcal F_k$, hence
\[
X(v)=av-Y(v)\in\mathcal F_k.
\]
It follows that $v\in\mathcal E_k$.

Given $v\in\mathcal E_k$ with $v_k=1$ and $a\in F$, consider the linear system
\begin{equation}\label{eq:linear-system}
Y(v)=av-X(v)
\end{equation}
in $Y\in\mathfrak b_r(F)$. It is easy to see that the map
\[
\mathfrak b_r(F)\longrightarrow\mathcal F_k,\qquad Y\longmapsto Y(v)
\]
is surjective. Hence the solution set of \eqref{eq:linear-system} is an affine space of dimension
\[
\dim\mathfrak b_r(F)-k.
\]
The set $\{v\in\mathcal E_k:v_k=1\}$ is either empty or an affine hyperplane of $\mathcal E_k$, of dimension $\dim\mathcal E_k-1$. It follows from the above discussion and \eqref{eq:Ek-dimension} that the incidence variety
\[
\mathcal V_k:=\{(Y,v,a)\in\mathfrak b_r(F)\times\mathcal E_k\times F:
v_k=1,\ Y(v)=av-X(v)\}
\]
has dimension at most
\[
(\dim\mathfrak b_r(F)-k)+(\dim\mathcal E_k-1)+1<\dim\mathfrak b_r(F).
\]
Let $\mathcal Y_k$ be the image of $\mathcal V_k$ in $\mathfrak b_r(F)$ under projection. Then $\bigcup_{k=1}^r\overline{\mathcal Y_k}$ is a proper closed subset of $\mathfrak b_r(F)$, and $\lambda$ is $(X+Y)$-cyclic for every
\[
Y\in\mathfrak b_r(F)\setminus\bigcup_{k=1}^r\overline{\mathcal Y_k}.
\]
This finishes the proof of the lemma.
\end{proof}

\section{Extension of characters}\label{sec:extension}
Let $N\subset\GL_n(F)$ be as in Theorem~\ref{thm:main}. Let $\mathbf U_n$ be the upper triangular maximal unipotent subgroup of $\GL_n$ with Lie algebra $\mathfrak u_n$, and let $U_n:=\mathbf U_n(F)$. By conjugation, from now on we assume that
\[
N\subset U_n.
\]
Since $p\ge n$ when $\operatorname{char}F=p$, in all cases we have the homeomorphisms
\[
\log:U_n\longrightarrow\mathfrak u_n(F)\quad\text{and}\quad
\exp:\mathfrak u_n(F)\longrightarrow U_n
\]
defined by the usual Taylor expansions. Put
\begin{equation}\label{eq:log-N}
\mathcal N:=\log(N)\subset\mathfrak u_n(F).
\end{equation}
The closed set $\mathcal N$ is an $F$-Lie subalgebra of $\mathfrak u_n(F)$ if $\operatorname{char}F=0$, and is an $\mathbb F_p$-Lie subalgebra of $\mathfrak u_n(F)$ if $\operatorname{char}F=p$ by the Lazard correspondence; see \cite[Chapitre~II, Theorem~4.6]{Laz54}.
Moreover, $\psi\circ\exp$ is an additive character of $\mathcal N$. 

\begin{remark}
Assume that $F$ is a local field, and $N=\mathbf N(F)$ with $\mathbf N$ a smooth connected subgroup of $\mathbf U_n$. If $\operatorname{char}F=0$, then
\begin{equation}\label{eq:Lie-N}
\mathcal N=\operatorname{Lie}(\mathbf N)(F).
\end{equation}
However if $\operatorname{char}F=p$, then \eqref{eq:Lie-N} is not true in general and $\mathcal N$ may not be an $F$-vector subspace of $\mathfrak u_n(F)$.
\end{remark}

Let $P_n\subset\GL_n(F)$ be the mirabolic subgroup consisting of matrices whose last row is $(0,\ldots,0,1)$. Then we have the Lie algebra
\[
\mathfrak p_n(F)=\left\{\begin{pmatrix}Y&v\\0&0\end{pmatrix}:
Y\in\mathfrak{gl}_{n-1}(F),\ v\in F^{n-1}\right\}.
\]
We identify its linear dual $\mathfrak p_n(F)^*=\mathfrak{gl}_{n-1}(F)\oplus F_{n-1}$ via the pairing
\begin{equation}\label{eq:pairing}
\left\langle(X,\lambda),\begin{pmatrix}Y&v\\0&0\end{pmatrix}\right\rangle
:=\operatorname{tr}(XY)+\lambda(v),
\quad X\in\mathfrak{gl}_{n-1}(F),\ \lambda\in F_{n-1}.
\end{equation}
For $k=1,2,\ldots,n-1$, let $N_k$ be the unipotent radical of the standard maximal parabolic subgroup of $\GL_n(F)$ of type $(k,n-k)$, with Lie algebra $\mathfrak n_k(F)$. Then it is easy to see that the annihilator of $\mathfrak n_k(F)$ in $\mathfrak p_n(F)^*$ is equal to
\begin{equation}\label{eq:annihilator}
\mathfrak n_k(F)^\perp=\{(X,\lambda)\in \mathfrak p_n(F)^*:X(\mathcal F_k)\subset\mathcal F_k,\ \lambda|_{\mathcal F_k}=0\},
\end{equation}
where $\mathcal F_k$ is given by \eqref{eq:standard-flag} with $r:=n-1$.

\begin{lemma}\label{lem:character-extension}
Let $E$ be a finite-dimensional $F$-vector space, let $L\subset E$
be a closed additive subgroup, and let $\eta:L\to\C^\times$ be a
smooth additive character. Put
\[
 \mathcal A_\eta=\{\phi\in E^*:\psi_F\circ\phi|_L=\eta\},
 \qquad
 \operatorname{Ann}_{\psi_F}(L)
 =\{\phi\in E^*:\psi_F(\phi(L))=1\}.
\]
Then $\mathcal A_\eta$ is nonempty and, for every
$\phi_0\in\mathcal A_\eta$,
\[
 \mathcal A_\eta=\phi_0+\operatorname{Ann}_{\psi_F}(L).
\]
For every $F$-vector subspace $D\subset E$, with its $F$-linear
annihilator $D^\perp\subset E^*$, one has
\begin{equation}\label{eq:character-annihilator}
 \mathcal A_\eta\subset D^\perp
 \quad\Longleftrightarrow\quad
 D\subset L\ \text{and}\ \eta|_D=1.
\end{equation}
If $L$ is an $F$-vector subspace, then $\mathcal A_\eta$ is an affine
space over $F$. In characteristic $p$ it is always an affine space
over $\mathbb F_p$.
\end{lemma}

\begin{proof}
Characters of closed subgroups of locally compact abelian groups
extend to the ambient group, and we have the isomorphism 
\[
 E^*\longrightarrow\widehat E,\qquad
 \phi\longmapsto\psi_F\circ\phi
\]
where $\widehat E$ is the Pontryagin dual of $E$. Thus $\mathcal A_\eta\ne\varnothing$,
and the formula for its translates follows. We use the character-extension
and double-annihilator theorems in \cite[Proposition~4.38 and
Corollary~4.41]{F95}.

Since $D$ is an $F$-vector subspace, $\psi_F\circ\phi$ is trivial on
$D$ if and only if $\phi|_D=0$. Consequently, if
$\mathcal A_\eta\subset D^\perp$, every element of the Pontryagin
annihilator of $L$ is trivial on $D$. The double-annihilator theorem
and the closedness of $L$ give $D\subset L$. Restriction of any
$\phi_0\in\mathcal A_\eta$ then gives $\eta|_D=1$. The converse
follows by restricting every extension to $D$.
If $L$ is $F$-linear, its Pontryagin annihilator is its $F$-linear
annihilator. In characteristic $p$ it is an $\mathbb F_p$-vector
subspace in all cases. This proves the last assertions.
\end{proof}

The main result of this section is as follows.
\begin{proposition}\label{prop:extension}
Let $F$ and $(N,\psi)$ be as in Theorem~\ref{thm:main} such that $N\subset U_n$.  Assume that there is no $k\in\{1,2,\ldots,n-1\}$ such that $N_k\subset N$ and $\psi|_{N_k}=1$. Then there exists a maximal unipotent subgroup $U\subset P_n$ and a nondegenerate character $\theta:U\to\C^\times$ such that
\[
\theta|_{U\cap N}=\psi|_{U\cap N}.
\]
\end{proposition}
\begin{proof}
The assertion is trivial when $n=1$. Assume that $n\ge2$. Put
\[
\mathcal A_\psi:=\{\phi\in\mathfrak p_n(F)^*:
\psi_F(\phi(x))=\psi(\exp(x))\text{ for all }x\in\mathcal N\}.
\]
By Lemma~\ref{lem:character-extension}, with
$E=\mathfrak p_n(F)$, $L=\mathcal N$ and $\eta=\psi\circ\exp$,
this is a nonempty affine space over $F$ in characteristic zero and
over $\mathbb F_p$ in characteristic $p$.
Since $\exp(\mathfrak n_k(F))=N_k$, equation
\eqref{eq:character-annihilator} gives, for $k=1,\ldots,n-1$,
\[
\mathcal A_\psi\subset\mathfrak n_k(F)^\perp
\quad\text{if and only if}\quad N_k\subset N\text{ and }\psi|_{N_k}=1.
\]
By assumption, each $\mathcal A_\psi\cap\mathfrak n_k(F)^\perp$
is empty or a proper affine subspace of $\mathcal A_\psi$.
In characteristic zero a finite union of proper affine subspaces
cannot cover an affine space over the infinite field $F$.
In characteristic $p$, choose a point outside each of these subspaces
and take their finite-dimensional affine span. Each intersection with
that span has at most $1/p$ of its points, and there are at most
$n-1<p$ intersections. Thus in all cases there exists an element
\begin{equation}\label{eq:choose-cyclic}
(X,\lambda)\in\mathcal A_\psi\setminus\bigcup_{k=1}^{n-1}\mathfrak n_k(F)^\perp.
\end{equation}
From \eqref{eq:pairing} it is clear that $(X+X',\lambda)\in\mathcal A_\psi$ for every $X'\in\mathfrak b_{n-1}(F)$. Hence by Lemma~\ref{lem:cyclic}, \eqref{eq:annihilator} and \eqref{eq:choose-cyclic}, we may assume that $(X,\lambda)\in\mathcal A_\psi$ and $\lambda$ is $X$-cyclic.

Define the following subspaces of $F^{n-1}$:
\[
\mathcal F'_k:=\bigcap_{i=0}^{n-k-2}\ker(\lambda X^i),
\quad k=0,1,\ldots,n-2;\qquad \mathcal F'_{n-1}:=F^{n-1}.
\]
Since $\lambda$ is $X$-cyclic, we have $\dim_F\mathcal F'_k=k$. Moreover,
\[
X(\mathcal F'_k)\subset\mathcal F'_{k+1}\quad\text{for }k=0,1,\ldots,n-2.
\]
For each $k=1,2,\ldots,n-2$, the induced map
\[
\mathcal F'_k/\mathcal F'_{k-1}\longrightarrow\mathcal F'_{k+1}/\mathcal F'_k
\]
is nonzero. Otherwise, $X(\mathcal F'_k)\subset\mathcal F'_k\subset\ker(\lambda)$, which implies that all the linear functionals $\lambda X^i$ vanish on the nonzero space $\mathcal F'_k$, contradicting that $\lambda$ is $X$-cyclic.

View $F^{n-1}$ as a subspace of $F^n$ via $v\mapsto(v,0)$. Then we have a complete flag
\[
\mathcal F'_\bullet:\{0\}\subset\mathcal F'_1\subset\cdots\subset
\mathcal F'_{n-1}\subset\mathcal F'_{n-1}+Fe_n=F^n,
\]
where $e_n$ is the last standard basis vector of $F^n$. Let $U$ be the maximal unipotent subgroup of $\GL_n(F)$ stabilizing $\mathcal F'_\bullet$. Clearly $U\subset P_n$, and with $\phi:=(X,\lambda)$ the map
\[
\theta:U\longrightarrow\C^\times,\qquad u\longmapsto\psi_F(\phi(\log u))
\]
is a nondegenerate character. By definition of $\mathcal A_\psi$, we have $\theta|_{U\cap N}=\psi|_{U\cap N}$. This proves the proposition.
\end{proof}

We conclude this section with a structural fact for general connected reductive groups. In the finite-field case, the inclusion-reversing property is proved by Grodal in the proof of \cite[Lemma~4.2, equation~(4.1)]{Gr02} and stated explicitly in \cite[Remark~4.3]{Gr02}, pp.~426--427. We give a relative-root proof which also applies over nonarchimedean local fields.

\begin{theorem}\label{thm:parabolic-radicals}
Let $F$ be a finite field or a nonarchimedean local field, let $\mathbf G$ be a connected reductive $F$-group, and let $\mathbf P,\mathbf Q$ be $F$-parabolic subgroups of $\mathbf G$. Put
\[
 \mathbf U_P=R_u(\mathbf P),\qquad \mathbf U_Q=R_u(\mathbf Q),
 \qquad U_P=\mathbf U_P(F),\quad U_Q=\mathbf U_Q(F).
\]
Then
\[
 \mathbf Q\subseteq\mathbf P
 \quad\Longleftrightarrow\quad
 \mathbf U_P\subseteq\mathbf U_Q
 \quad\Longleftrightarrow\quad
 U_P\subseteq U_Q.
\]
In particular, the unipotent radical determines its parabolic subgroup, and strict inclusions of unipotent radicals correspond to strict reverse inclusions of parabolic subgroups.
\end{theorem}
\begin{proof}
By the relative Bruhat decomposition, $\mathbf P$ and $\mathbf Q$ contain a common maximal $F$-split torus $\mathbf A$. Put $\mathbf M=Z_{\mathbf G}(\mathbf A)$. Let $\Phi_{\mathrm{nd}}$ be the set of indivisible relative roots, and for $\alpha\in\Phi_{\mathrm{nd}}$ let $\mathbf U_\alpha$ be the relative root group, including the $2\alpha$-part when it occurs. For an $F$-parabolic subgroup $\mathbf R$ containing $\mathbf A$, set
\[
 \Psi_R=\{\alpha\in\Phi_{\mathrm{nd}}:\mathbf U_\alpha\subseteq\mathbf R\},
 \qquad
 \Theta_R=\{\alpha\in\Phi_{\mathrm{nd}}:\mathbf U_\alpha\subseteq R_u(\mathbf R)\}.
\]
The relative root description of parabolic subgroups \cite[\S\S4--5]{BT65} gives
\[
 \Psi_R=\Phi_{\mathrm{nd}}\setminus(-\Theta_R).
\]
Moreover, $\mathbf R$ is generated by $\mathbf M$ and the root groups with roots in $\Psi_R$, and $R_u(\mathbf R)$ is generated by those with roots in $\Theta_R$. If $\alpha\notin\Theta_R$, the intersection $\mathbf U_\alpha\cap R_u(\mathbf R)$ is trivial.

Suppose $U_P\subseteq U_Q$. Each $\mathbf U_\alpha$ is a nontrivial $F$-split unipotent group, so $\mathbf U_\alpha(F)\ne\{1\}$. For $\alpha\in\Theta_P$, we have $\mathbf U_\alpha(F)\subseteq U_P\subseteq U_Q$; the intersection assertion forces $\alpha\in\Theta_Q$. Thus $\Theta_P\subseteq\Theta_Q$, and taking complements of their negatives gives $\Psi_Q\subseteq\Psi_P$. The generation statement now implies $\mathbf Q\subseteq\mathbf P$.

Conversely, $\mathbf Q\subseteq\mathbf P$ implies $\Psi_Q\subseteq\Psi_P$, hence $\Theta_P\subseteq\Theta_Q$ and therefore $\mathbf U_P\subseteq\mathbf U_Q$. This implies the corresponding inclusion on $F$-points. The final assertions follow from the equivalences.
\end{proof}

The following corollary of Theorem~\ref{thm:parabolic-radicals} explains why in Proposition~\ref{prop:extension} it suffices to single out the subgroups $N_k$ contained in $N$, rather than consider all the unipotent radicals of parabolic subgroups of $\GL_n$ contained in $N$, as suggested by Conjecture~\ref{conj:degenerate}.

\begin{corollary}\label{cor:standard-radicals}
Let $F$ be a finite field or a nonarchimedean local field, and let $B_n$ and $U_n$ be the upper triangular Borel subgroup of $\GL_n(F)$ and its unipotent radical. If $\mathbf P$ is an $F$-parabolic subgroup of $\GL_n$ such that $R_u(\mathbf P)(F)\subseteq U_n$, then $\mathbf P$ is standard. Thus the unipotent radicals of proper parabolic subgroups of $\GL_n$ contained in $U_n$ are exactly the unipotent radicals of proper standard parabolic subgroups.

In particular, if a subgroup $H\subseteq U_n$ contains $R_u(\mathbf P)(F)$ for some proper $F$-parabolic subgroup $\mathbf P$, then it contains $N_k$ for some $1\leq k\leq n-1$, where
\[
 N_k=\left\{\begin{pmatrix}1_k&X\\0&1_{n-k}\end{pmatrix}:X\in\operatorname{Mat}_{k\times(n-k)}(F)\right\}
\]
is the unipotent radical of the standard maximal parabolic subgroup of type $(k,n-k)$.
\end{corollary}
\begin{proof}
Theorem~\ref{thm:parabolic-radicals}, applied to $\mathbf P$ and the upper triangular Borel subgroup, gives $\mathbf B_n\subseteq\mathbf P$. If $\mathbf P$ is proper, choose a standard maximal parabolic subgroup $\mathbf P_k$ containing it. The same theorem gives $N_k=R_u(\mathbf P_k)(F)\subseteq R_u(\mathbf P)(F)\subseteq H$.
\end{proof}

\section{Proof of Theorem \ref{thm:main}}\label{sec:main-proof}
We will apply the bottom piece of the Bernstein--Zelevinsky filtration to prove Theorem~\ref{thm:main}. Let $U$ be a maximal unipotent subgroup of $P_n$, and $\theta:U\to\C^\times$ be a nondegenerate character. By \cite[\S5.15]{BZ76} and \cite[Proposition~1.2]{Gel70}, for every irreducible generic representation $\pi$ of $\GL_n(F)$, we have a $P_n$-equivariant embedding
\begin{equation}\label{eq:bottom-piece}
\cind_U^{P_n}\theta\hookrightarrow\pi|_{P_n},
\end{equation}
where $\cind$ denotes compact induction.

\begin{lemma}\label{lem:compatible-character}
Let $(N,\psi)$ be as in Theorem~\ref{thm:main} such that $N\subset U_n$. If $\theta|_{U\cap N}=\psi|_{U\cap N}$, then
\[
\bigl(\cind_U^{P_n}\theta\bigr)_{N,\psi}\ne\{0\}.
\]
\end{lemma}
\begin{proof}
The assertion is immediate when $F$ is finite. Assume that $F$ is local. Apply Lemma~\ref{lem:closed-orbit} below to the action of $U\times N$ on the affine group $P_n$ by $(u,n)g=ugn^{-1}$, at $g=1$. Then $UN$ is closed in $P_n$, and the orbit homeomorphism, followed by passage to left $U$-cosets, shows that the natural map
\[
U\cap N\backslash N\longrightarrow U\backslash UN\subset U\backslash P_n
\]
is a homeomorphism onto a closed subspace of $U\backslash P_n$. It follows that the map
\[
\cind_U^{P_n}\theta\longrightarrow\C,\qquad
f\longmapsto\int_{U\cap N\backslash N}f(x)\overline{\psi}(x)\,dx
\]
is well-defined and is a nonzero $(N,\psi)$-equivariant continuous linear functional, where $dx$ is a fixed right $N$-invariant measure and the integral converges absolutely.
\end{proof}

To complete the proof of Lemma~\ref{lem:compatible-character}, we needed to show that the product $UN$ is a closed subset of $P_n$. This does not follow from general facts about topological groups, even when $U$ and $N$ are closed. For this purpose, we use the following lemma on the closedness of rational unipotent orbits in affine varieties, including in positive characteristic, for which we sketch a proof. Note that Bernstein--Zelevinsky \cite[\S6.15, Theorem~A, and the Appendix, A.2]{BZ76} give constructibility of rational orbits and hence local closedness (see also \cite[\S6.8(c)]{BZ76}). 

\begin{lemma}\label{lem:closed-orbit}
Let $F$ be a nonarchimedean local field. Let $\mathbf H$ be a smooth connected unipotent $F$-algebraic group, acting on an affine $F$-variety $\mathbf X$. Let $x\in\mathbf X(F)$, and let $\mathbf S$ be the scheme-theoretic stabilizer of $x$ in $\mathbf H$. Put $H:=\mathbf H(F)$ and $S:=\mathbf S(F)$. Then $Hx$ is closed in $\mathbf X(F)$, and the natural map
\[
H/S\longrightarrow Hx
\]
is a homeomorphism.
\end{lemma}
\begin{proof}
The geometric orbit of a unipotent group is closed by the Kostant--Rosenlicht theorem (see \cite[Proposition~4.10]{B91}). Let $\mathbf O=\mathbf H/\mathbf S$ be the orbit of $x$ with its homogeneous-space structure. Then the morphism $\mathbf H\to\mathbf O$ is an fppf $\mathbf S$-torsor. The scheme-theoretic subgroup $\mathbf S$ of $\mathbf H$ has
reductive rank zero, whether or not it is smooth. By
\cite[Theorem~1.2(1)(c)]{GGM14}, the image of $\mathbf H(F)$ is closed
in $\mathbf O(F)$. Part~(2) of the same theorem identifies this image,
with its subspace topology, with $H/S$. Since $\mathbf O$ is a variety,
part~(3) also makes $H\to Hx$ a topological principal $S$-bundle.
Finally, $\mathbf O(F)$ is closed in $\mathbf X(F)$, which proves the lemma.
\end{proof}

\begin{proposition}\label{prop:nonvanishing}
Let $F$ and $(N,\psi)$ be as in Theorem~\ref{thm:main} such that $N\subset U_n$. If $\ker\psi$ contains none of $N_1,N_2,\ldots,N_{n-1}$, then
\[
\pi_{N,\psi}\ne\{0\}
\]
for every irreducible generic representation $\pi$ of $\GL_n(F)$.
\end{proposition}
\begin{proof}
By Proposition~\ref{prop:extension}, there exist a maximal unipotent subgroup $U$ of $P_n$ and a nondegenerate character $\theta$ of $U$ such that $\theta|_{U\cap N}=\psi|_{U\cap N}$. The proposition follows from \eqref{eq:bottom-piece}, Lemma~\ref{lem:compatible-character} and the exactness of the functor $\sigma\mapsto\sigma_{N,\psi}$ on the category of smooth representations of $N$.
\end{proof}

We need the following variant with similar proof of Proposition~\ref{prop:nonvanishing}.

\begin{proposition}\label{prop:product}
Let $F$ be as in Theorem~\ref{thm:main}. Assume that $d_1,\ldots,d_r$ are positive integers with $\sum_{j=1}^r d_j=n$. Put
\[
M=\prod_{j=1}^r\GL_{d_j}(F),\qquad U_M=\prod_{j=1}^r U_{d_j}.
\]
For $j=1,\ldots,r$ and $k=1,\ldots,d_j-1$, let $N_{j,k}$ be the unipotent radical of the standard maximal parabolic subgroup of $\GL_{d_j}(F)$ of type $(k,d_j-k)$, embedded into the $j$-th factor of $M$. Assume that $N_M$ is a subgroup of $U_M$, and is of the form $\mathbf N_M(F)$ for a smooth connected unipotent $F$-algebraic group $\mathbf N_M$ when $F$ is local. Let $\psi_M$ be a character of $N_M$ such that $\ker\psi_M$ contains none of the groups $N_{j,k}$. Then
\[
\sigma_{N_M,\psi_M}\ne\{0\}
\]
for every irreducible generic representation $\sigma$ of $M$.
\end{proposition}

\begin{proof}
The subgroup $N_M$ need not be a product. We therefore choose one
functional extending its character before choosing flags in the
individual blocks. Put
\[
 P:=\prod_{j=1}^r P_{d_j},\qquad
 \mathfrak p(F):=\bigoplus_{j=1}^r\mathfrak p_{d_j}(F),\qquad
 \mathcal N_M:=\log N_M\subset \mathfrak p(F),
\]
where logarithms are taken block by block.
The logarithmic argument of Section~\ref{sec:extension}, applied
block by block, shows that $\mathcal N_M$ is a closed additive subgroup
of $\mathfrak p(F)$ and that $\psi_M\circ\exp$ is an additive character.
Write
\[
 \mathcal A_{\psi_M}:=\{\phi\in \mathfrak p(F)^*:\psi_F(\phi(x))=\psi_M(\exp x)
                           \text{ for every }x\in  \mathcal N_M\},
\]
where $\mathfrak p(F)^*$ denotes the $F$-linear dual.
Lemma~\ref{lem:character-extension} shows that $\mathcal A_{\psi_M}$
is nonempty and affine over $F$ in characteristic zero and over
$\mathbb F_p$ in characteristic $p$.

For $j=1,\ldots, r$ and $k=1,\ldots, d_j-1$, let
$\mathfrak n_{j,k}(F)$ be the Lie algebra of $N_{j,k}$, embedded into the $j$-th summand of  $\mathfrak p(F)$, and let $\mathfrak n_{j,k}(F)^{\perp}\subset \mathfrak p(F)^*$ be its
annihilator. Equation~\eqref{eq:character-annihilator} and
$N_{j,k}\not\subset\ker\psi_M$ give
$\mathcal A_{\psi_M}\not\subset\mathfrak n_{j,k}(F)^\perp$.
There are $\sum_j(d_j-1)=n-r$ such subspaces. The affine-avoidance
argument in Proposition~\ref{prop:extension}, using $n-r<p$ in
characteristic $p$, therefore gives 
\[
\phi = (\phi_j)_j \in \mathcal A_{\psi_M} \setminus \bigcup_{j=1,\ldots, r, \, k=1,\ldots, d_j-1} \mathfrak n_{j,k}(F)^\perp.
\]

For a block with $d_j=1$, take $U'_j=\{1\}$ and $\theta_j=1$.
For each remaining block, write $\phi_j=(X_j,\lambda_j)$ under
\eqref{eq:pairing}. By \eqref{eq:annihilator} and Lemma~\ref{lem:cyclic},
we may add an upper triangular matrix to $X_j$ so that $\lambda_j$
is $X_j$-cyclic. These changes leave $\phi|_{\mathcal N_M}$ unchanged: the trace
pairing of an upper triangular matrix with a strictly upper triangular
matrix is zero. The flag construction in the proof of
Proposition~\ref{prop:extension} now gives a maximal unipotent
subgroup $U'_j\subset P_{d_j}$ on which
$\theta_j(u)=\psi_F(\phi_j(\log u))$ is a nondegenerate character.
 Put
$U':=\prod_j U'_j$ and $\theta:=\boxtimes_j\theta_j$. For
$u=(u_j)_j\in U'\cap N_M$ we have
\[
 \theta(u)=\psi_F\!\left(\sum_j\phi_j(\log u_j)\right)
          =\psi_F(\phi(\log u))=\psi_M(u).
\]

Write $\sigma=\boxtimes_j\sigma_j$, with $\sigma_j$ irreducible
generic on $\GL_{d_j}(F)$. Taking the external tensor product of
the embeddings \eqref{eq:bottom-piece} gives
\[
 \cind_{U'}^P\theta
 \simeq\boxtimes_j\cind_{U'_j}^{P_{d_j}}\theta_j
 \hookrightarrow \sigma|_P .
\]
The proof of Lemma~\ref{lem:compatible-character} applies to this
product: integration of $f(x)\overline\psi_M(x)$ over
$(U'\cap N_M)\backslash N_M$ defines a nonzero
$(N_M,\psi_M)$-equivariant functional on $\cind_{U'}^P\theta$.
Exactness of $(N_M,\psi_M)$-coinvariants
therefore embeds the nonzero module $(\cind_{U'}^{P}\theta)_{N_M,\psi_M}$
in $\sigma_{N_M,\psi_M}$.
\end{proof}

For every subset
\[
I=\{i_1<i_2<\cdots<i_{r-1}\}\subset\{1,2,\ldots,n-1\},
\]
set $i_0:=0$, $i_r:=n$ and let $Q_I$ be the standard parabolic subgroup of $\GL_n(F)$ of type $(d_1,d_2,\ldots,d_r)$, where
\[
d_j:=i_j-i_{j-1},\qquad j=1,2,\ldots,r.
\]
Let $U_I$ be the unipotent radical of $Q_I$, and $M_I=\prod_{j=1}^r\GL_{d_j}(F)$ be the block diagonal Levi subgroup of $Q_I$. Recall that $N_i$, $i=1,2,\ldots,n-1$, denotes the unipotent radical of the standard maximal parabolic subgroup of type $(i,n-i)$. Then $U_I$ is generated by $N_i$, $i\in I$.

For a subgroup $N\subset U_n$ and a character $\psi$ of $N$, put
\begin{equation}\label{eq:I-N-psi}
I(N,\psi):=\{1\leq i<n : N_i\subset\ker\psi\}.
\end{equation}
Now we can prove the following result, which implies Theorem~\ref{thm:main}.
\begin{theorem}\label{thm:precise}
Let $F$ and $(N,\psi)$ be as in Theorem~\ref{thm:main} such that $N\subset U_n$, and let $\pi$ be an irreducible generic representation of $\GL_n(F)$. Then with $I(N,\psi)$ given by \eqref{eq:I-N-psi},
\[
\pi_{N,\psi}=\{0\}\quad\text{if and only if}\quad \pi_{U_{I(N,\psi)}}=\{0\}.
\]
\end{theorem}
\begin{proof}
Write $U:=U_{I(N,\psi)}$, $M:=M_{I(N,\psi)}$ and $Q:=Q_{I(N,\psi)}$. Since $U$ is generated by $N_i$ with $i\in I(N,\psi)$ so that $N_i\subset\ker\psi$, we have $U\subset\ker\psi$. The Levi decomposition $Q=M\ltimes U$ restricts to a semi-direct product decomposition
\[
N\simeq(N\cap M)\ltimes U.
\]
Write $N_M:=N\cap M$ and $\psi_M:=\psi|_{N_M}$.
When $F$ is local, the inclusion $U\subset N$ on $F$-points gives
$\mathbf U\subset\mathbf N$ because $\mathbf U(F)$ is Zariski dense
in $\mathbf U$. The Levi projection identifies
$\mathbf N\cap\mathbf M$ with the smooth connected quotient
$\mathbf N/\mathbf U$. Thus $N_M=(\mathbf N\cap\mathbf M)(F)$
satisfies the hypotheses of Proposition~\ref{prop:product}.
Then
\[
\pi_{N,\psi}\simeq(\pi_U)_{N_M,\psi_M}.
\]
Thus if $\pi_U=\{0\}$, then $\pi_{N,\psi}=\{0\}$.

Conversely, assume that $\pi_U\ne\{0\}$, which is a smooth representation of $M$ of finite length.  Let $\sigma$ be an irreducible quotient of $\pi_U$.
By Frobenius reciprocity,
\[
 \operatorname{Hom}_M(\pi_U,\sigma)
 \simeq\operatorname{Hom}_G(\pi,\Ind_{Q}^G\sigma).
\]
Thus $\pi$ embeds into $\Ind_{Q}^G\sigma$.
Choose a nondegenerate character $\theta$ of $U_n$ for which
$\pi_{U_n,\theta}\ne0$. Exactness of twisted coinvariants gives
\[
 0\ne\pi_{U_n,\theta}
 \hookrightarrow
 (\Ind_{Q}^G\sigma)_{U_n,\theta}.
\]
Heredity of Whittaker models under parabolic induction now implies
that $\sigma$ is generic. For local fields this is \cite[\S V, Theorem~4]{R75}. For finite
fields it follows from
Harish--Chandra restriction of the Gelfand--Graev representation,
\cite[Theorem~2.1]{BR06}.

Write $I(N,\psi)=\{i_1<\cdots<i_{r-1}\}$ and $i_0=0$, $i_r=n$.
If $N_{j,k}\subset\ker\psi_M$, put $i=i_{j-1}+k$.
The Levi projection maps $N_i$ onto $N_{j,k}$ and has kernel
$N_i\cap U$. Since both $U$ and $N_{j,k}$ are contained in
$\ker\psi$, this gives $N_i\subset\ker\psi$, contradicting
$i\notin I(N,\psi)$. Proposition~\ref{prop:product} therefore gives
$\sigma_{N_M,\psi_M}\ne0$. The quotient $\pi_U\twoheadrightarrow\sigma$
induces a surjection
\[
 \pi_{N,\psi}\simeq(\pi_U)_{N_M,\psi_M}
 \twoheadrightarrow\sigma_{N_M,\psi_M},
\]
so $\pi_{N,\psi}\ne0$.

The theorem is proved.
\end{proof}

\section{Generalized Whittaker models}\label{sec:generalized}
In this section, let $F$ be any finite or nonarchimedean local field. 
Let $f\in\mathfrak g=\mathfrak{gl}_n(F)$ be nilpotent.
On each Jordan chain $v_1,\ldots,v_d$, where $fv_i=v_{i+1}$ for
$i<d$ and $fv_d=0$, assign the successive integer weights
$d-1,d-3,\ldots,1-d$. Write $V_b\subset F^n$ for the span of the
basis vectors of weight $b$. For $i\in\mathbb Z$ put
\[
 \mathfrak g_i=\{X\in\mathfrak g:X(V_b)\subseteq V_{b+i}
          \text{ for all }b\},\qquad
 \mathfrak g_{\ge i}=\bigoplus_{j\ge i}\mathfrak g_j .
\]
The alternating form $\omega_f(X,Y)=\operatorname{tr}(f[X,Y])$
on $\mathfrak g_1$ is nondegenerate in every characteristic.
Indeed, the trace pairing between $\mathfrak g_{-1}$ and $\mathfrak g_1$
is perfect, so the radical of $\omega_f$ is the kernel of $X\mapsto[f,X]$ on
$\mathfrak g_1$. A map commuting with $f$ from a chain of length $d$
to one of length $e$ sends the first vector into $\ker f^d$ on the
target. The largest weight in this kernel is $e-1-2\max(0,e-d)$,
which is at most $d-1$, the weight of the source vector. Such a map
of degree one must therefore vanish.

Put $G_{\geq 1}=1+\mathfrak g_{\ge1}$ and $G_{\geq 2}=1+\mathfrak g_{\ge2}$. For
$x\in\mathfrak g_{\ge1}$, write $x_i$ for its degree-$i$ component.
Then $G_{\geq 2}$ is normal in $G_{\geq 1}$ and has the $G_{\geq 1}$-invariant character
\[
 \chi_f(1+z)=\psi_F(\operatorname{tr}(fz_2)),\qquad z\in\mathfrak g_{\ge2}.
\]
For a Lagrangian $\mathfrak l\subseteq\mathfrak g_1$ put
$N_{f,\mathfrak l}=1+(\mathfrak l\oplus\mathfrak g_{\ge2})$,
the inverse image of $\mathfrak l$ under $G_{\geq 1}\to\mathfrak g_1$,
$1+x\mapsto x_1$. Let $\psi_{f,\mathfrak l}$ be any smooth character
of $N_{f,\mathfrak l}$ extending $\chi_f$. Such extensions exist:
for $x,y\in\mathfrak g_{\ge1}$ the group commutator of $1+x$ and $1+y$
is congruent to $1+[x_1,y_1]$ modulo $1+\mathfrak g_{\ge3}$.
By isotropy, $N_{f,\mathfrak l}/\ker\chi_f$ is therefore locally compact
abelian, and the character of its closed subgroup $G_{\geq 2}/\ker\chi_f$
extends by the character-extension theorem for locally compact abelian
groups \cite[Corollary~4.41]{F95}. The extension is smooth: a sufficiently
small compact open subgroup maps into an arc of the unit circle
containing no nontrivial subgroup, hence lies in its kernel.
For a smooth representation $\pi$ of $\GL_n(F)$, its generalized Whittaker model
associated to $f$ and a choice of the Lagrangian $\mathfrak l \subset \mathfrak g_1$
is the space 
\[
\operatorname{Hom}_{N_{f,\mathfrak l}}(\pi,\psi_{f,\mathfrak l}).
\]
For the finite-field construction by extension from degree two see
\cite[Lemma~5.13, Definition~5.14 and Lemma~5.15]{T16}.

When $2$ is invertible in $F$, one can take
\[
 \psi_{f,\mathfrak l}(1+x)=
 \psi_F\!\left(\operatorname{tr}\bigl(f(x_2-\tfrac12 x_1^2)\bigr)\right).
\]
Isotropy makes this a character, and where the logarithm is defined
the formula is precisely $\psi_F(\operatorname{tr}(f\log(1+x)))$, since
only degree two contributes to the trace. Thus our definition agrees
there with the usual polarized generalized Whittaker model; compare
\cite[Definition~2.5.1]{GGS17} in characteristic zero.

\begin{corollary}\label{cor:generalized}
Let $F$ be any finite or nonarchimedean local field and let $\pi$ be an
irreducible generic representation of $\GL_n(F)$. For every nilpotent
$f\in\mathfrak{gl}_n(F)$, every Lagrangian $\mathfrak l\subseteq\mathfrak g_1$,
and every smooth extension $\psi_{f,\mathfrak l}$ of $\chi_f$ to
$N_{f,\mathfrak l}$, one has
\[
 \pi_{N_{f,\mathfrak l},\psi_{f,\mathfrak l}}\ne\{0\}.
\]
\end{corollary}

When $F$ is local of characteristic zero, Corollary \ref{cor:generalized} also
follows from \cite[Theorem~E]{GGS17}. We give the following direct proof, which applies
in every characteristic.

\begin{proof}
List the Jordan chains consecutively, each in the order
$v_1,\ldots,v_d$, and put a longest chain last. Let $U'$ be the upper
unitriangular subgroup for this ordered basis and $P$ the corresponding
mirabolic subgroup. The last vector has
minimal weight among all the chains, so an endomorphism of positive
degree cannot have any component in its row. Thus $G_{\geq 1}\subset P$.
On $U'$ take the generic character
\[
 \theta(u)=\psi_F\!\left(\sum_{i=1}^{n-1}u_{i,i+1}\right).
\]
For $u\in U'\cap G_{\geq 2}$, the entries in this sum internal to the chains
give precisely $\operatorname{tr}(f(u-1))$. The remaining entries
join the last vector of one chain of length $d$ to the first of the
next, of length $e$: the corresponding matrix unit has degree
$(1-d)-(e-1)\le0$, so its coefficient in $u-1$ vanishes. Hence
$\theta=\chi_f$ on $U'\cap G_{\geq 2}$.

Write $N=N_{f,\mathfrak l}$, $\tau=\psi_{f,\mathfrak l}$, and
$K=N\cap U'$. We first arrange agreement on all of $K$. Since
$\theta/\tau$ is trivial on $K\cap G_{\geq 2}$, it descends under $1+x\mapsto x_1$
to a smooth additive character of
$W=\mathfrak l\cap\mathfrak u'$, where $\mathfrak u'$ is the space
of strictly upper triangular matrices in our ordered basis. Indeed, taking a
degree component preserves strict upper triangularity in the ordered
Jordan basis, and $X\mapsto 1+X$ gives a section on $W$.
By additive self-duality this character is $\psi_F\circ\lambda$ for
some $\lambda\in W^*$. Extend $\lambda$ $F$-linearly to $\mathfrak g_1$.
By nondegeneracy there is $Y\in\mathfrak g_1$ such that
$\lambda(X)=\omega_f(Y,X)$ for all $X\in\mathfrak g_1$.
The element $h=1+Y$ normalizes $N$, and the commutator calculation above
gives, for $k=1+x\in K$,
\[
 \tau(hkh^{-1})=\tau(k)\,\psi_F(\omega_f(Y,x_1))=\theta(k).
\]
Thus the conjugate character $\tau^h(u):=\tau(huh^{-1})$ agrees
with $\theta$ on $N\cap U'$.

The mirabolic embedding \eqref{eq:bottom-piece} gives
$\cind_{U'}^P\theta\hookrightarrow\pi|_P$ in every characteristic.
The argument of Lemma~\ref{lem:compatible-character} now gives a
nonzero $(N,\tau^h)$-functional on $\cind_{U'}^P\theta$ by integration
over $K\backslash N$. In the local case the required closed orbit
and quotient topology follow from Lemma~\ref{lem:closed-orbit} for
the action of $U'\times N$ on $P$; in the finite case the integral is
a sum. Exactness of twisted coinvariants gives $\pi_{N,\tau^h}\ne0$.
Since $h$ normalizes $N$, the action of $\pi(h)$ identifies these
coinvariants with $\pi_{N,\tau}$.
\end{proof}

The example of $\GL_3(\mathbb F_2)$ illustrates which characters of
$N_{f,\mathfrak l}$ extend the prescribed character $\chi_f$ of $G_{\geq 2}$.
Take $f=E_{31}$ with weights $(1,0,-1)$ on the standard basis,
and $\mathfrak l=\mathbb F_2 X$ with $X=E_{12}+E_{23}$, where $E_{ij}$
denotes the elementary matrix. Since $(1+X)^2=1+E_{13}$ and
$\chi_f(1+E_{13})=-1$, any such character satisfies
$\psi_{f,\mathfrak l}(1+X)^2=-1$ and hence takes the value $i$ or $-i$
on $1+X$. The character in Example~\ref{ex:characteristic-bound} takes
the value $-1$ on $1+X$, so its restriction to $G_{\geq 2}$ differs from
$\chi_f$.

\section{Restriction of the Steinberg representation}\label{sec:steinberg}

We recall Borel--Serre's resolution of the Steinberg representation for a connected reductive group \cite[\S3]{BS76}, and its realization as the space of compactly supported locally constant functions on the open Bruhat cell. The restriction to a minimal parabolic subgroup then gives the nonvanishing of twisted coinvariants in Corollary~\ref{cor:steinberg-conjecture}.

Let $F$ be a nonarchimedean local field or a finite field, and let $\mathbf G$ be an arbitrary connected reductive $F$-group, with $G=\mathbf G(F)$. Choose a maximal $F$-split torus $\mathbf A$, put $\mathbf M=Z_{\mathbf G}(\mathbf A)$, and choose opposite minimal $F$-parabolic subgroups $\mathbf B=\mathbf M\mathbf U$ and $\overline{\mathbf B}=\mathbf M\overline{\mathbf U}$. Write $B=MU$, and similarly for the other groups of $F$-points. Here $B$ denotes a minimal $F$-parabolic subgroup; it is a Borel subgroup only when $\mathbf G$ is quasi-split. The Levi group $M$ need not be a torus. Its derived group is $F$-anisotropic.

Let $\Phi=\Phi(\mathbf G,\mathbf A)$ be the relative root system, possibly nonreduced, with positive system determined by $U$. Its reduced system of indivisible roots determines the finite Coxeter system $(W,S)$, where
\[
 W=N_{\mathbf G}(\mathbf A)(F)/M.
\]
Write $\ell$ for Coxeter length, $\leq$ for Bruhat order, and $w_0$ for the longest element of $W$. Choose representatives of the elements of $W$ in $N_{\mathbf G}(\mathbf A)(F)$ and, by abuse of notation, denote the representative of $w$ again by $w$. The rational Bruhat decomposition and relative root subgroup factorizations used here are those of Borel--Tits \cite{BT65}; see also \cite[\S1]{C95} for the local-field setup.

For $I\subseteq S$, let $W_I=\langle I\rangle$ and let $P_I$ be the standard parabolic subgroup containing $B$. Thus $P_\varnothing=B$ and $P_S=G$. Set
\[
 X=B\backslash G,\qquad X_I=P_I\backslash G,\qquad q_I:X\longrightarrow X_I,\qquad E_I=C^\infty(X_I).
\]
Here $C^\infty$ means locally constant complex-valued functions. These flag spaces are compact, Hausdorff and totally disconnected; $\Cc$ means locally constant functions of compact support on a possibly noncompact stratum. Pullback embeds $E_I$ in $E_\varnothing$. We always use unnormalized induction and the right action
\[
 (\rho(g)f)(Hx)=f(Hxg).
\]
Define
\[
 R=\sum_{s\in S}q_{\{s\}}^*E_{\{s\}},\qquad \St_G=E_\varnothing/R.
\]
\subsection{The augmented parabolic chain complex}
Fix a total order on $S$ and put $r=|S|$. For $0\leq i\leq r$ let
\[
 K_i=\bigoplus_{\substack{I\subseteq S\\|I|=i}}E_I .
\]
Thus $K_r=E_S=\C$ and $K_0=E_\varnothing=C^\infty(B\backslash G)$.
For $J\subseteq I$ write $p_{J,I}:X_J\to X_I$ for the projection
$P_Jg\mapsto P_Ig$. Write $f[I]$ for $f\in E_I$ in its summand of
$K_i$. If $I=\{s_1<\cdots<s_i\}$ and $i\geq1$, the boundary is
\[
 d_i(f[I])=
 \sum_{j=1}^i(-1)^{j-1}
 \bigl(p_{I\setminus\{s_j\},I}^*f\bigr)[I\setminus\{s_j\}].
\]
Deleting two elements in either order gives the same pullback with
opposite signs, so $d_{i-1}d_i=0$ for $i\geq2$. The image of
$d_1((f_s)_{s\in S})=\sum_{s\in S}q_{\{s\}}^*f_s$ is $R$ when $r\geq1$.
With $\varepsilon:K_0\to\St_G$ the defining quotient, the augmented
complex of smooth $G$-modules is
\begin{equation}\label{eq:restriction-2-1}
 0\longrightarrow K_r\xrightarrow{d_r}K_{r-1}
 \longrightarrow\cdots\longrightarrow K_1\xrightarrow{d_1}K_0
 \xrightarrow{\varepsilon}\St_G\longrightarrow0 .
\end{equation}
For $r=0$ this reads $0\to\C\xrightarrow{\mathrm{id}}\C\to0$.

\begin{theorem}[Borel--Serre]\label{thm:steinberg}\leavevmode
\begin{enumerate}
\item \textup{Exactness.} The augmented parabolic chain complex
\eqref{eq:restriction-2-1} is exact.
\item \textup{Open cell.} Let $\Omega=B\backslash Bw_0B$. Extension by
zero to $B\backslash G$, followed by the quotient map $\varepsilon$, gives
a $B$-equivariant isomorphism
\[
 \Cc(\Omega)\xrightarrow{\ \sim\ }\St_G|_B.
\]
Consequently $\St_G|_B\simeq\cind_M^B\one$ and $\St_G|_U\simeq\Cc(U)$,
with the right regular action of $U$.
\end{enumerate}
\end{theorem}

\begin{proof}
For the sake of completeness, we outline the proof of Borel--Serre.
Both the exactness of the augmented complex and the contribution of the
open cell follow from the same filtration by Bruhat cells
\cite[Corollary~3.3, Proposition~3.5 and Remark~3.6]{BS76}.
Borel--Serre work over a nonarchimedean local field.
Passage to the adjoint group identifies the flag varieties of every type
and their projections $G$-equivariantly, giving the connected reductive
case. The same filtration works over a finite field.

Consider $K_\bullet$ before augmentation. Order the Bruhat cells
$C_w=B\backslash BwB$ in reverse Bruhat order, beginning with $C_{w_0}=\Omega$.
Filter each $E_I$ by requiring the support of its pullback to $B\backslash G$
to lie in the successive open unions of these cells. The boundary maps
preserve the filtration because pulling back to $B\backslash G$ gives the
same function before and after a parabolic pullback.

To compute the quotient at $C_w$, consider its image in $P_I\backslash G$.
The Bruhat cells in $B\backslash G$ above this $B$-orbit are precisely the
$C_v$ with $v\in W_Iw$. The cell indexed by the minimal representative
of $W_Iw$ enters the filtration last and maps isomorphically onto the orbit. A pullback can be nonzero above
the orbit only once the whole fibre is allowed in its support. Restriction
to the new stratum then gives all its compactly supported locally constant
functions \cite[\S5.11]{BZ77}. Thus $E_I$ contributes at $w$ precisely
when $w$ is minimal in $W_Iw$, or equivalently when
$I\subseteq A(w)$ \cite[Ch.~IV, \S1]{Bo02}. Here
\[
 A(w)=\{s\in S:\ell(sw)>\ell(w)\}
\]
is the set of simple reflections $s$ for which multiplying $w$ on the left
by $s$ increases its length. In degree $i$ this gives
\[
 \gr_w K_i\simeq\bigoplus_{\substack{I\subseteq A(w)\\ |I|=i}}\Cc(C_w)[I].
\]
Each boundary pullback acts as the identity on $\Cc(C_w)$. Its sign is
the sign for deleting a vertex, so this is the augmented chain complex
of the simplex with vertex set $A(w)$, with coefficients in
$\Cc(C_w)$ and the empty face in degree zero.

If $w\ne w_0$, a reduced expression for $w$ can be extended on the left
to one for $w_0$. Thus $A(w)$ is nonempty and its simplex is contractible,
so the graded complex is acyclic.
Since $w_0$ has maximal length, $A(w_0)=\varnothing$, so the first filtration term consists
of $\Cc(\Omega)$ in degree zero alone. Each subsequent quotient is acyclic.
The long exact homology sequence at each step therefore shows that inclusion
of the first term induces
\[
 H_i(K_\bullet)=0\quad(i>0),\qquad
 \Cc(\Omega)\xrightarrow{\ \sim\ }H_0(K_\bullet).
\]
Since $H_0(K_\bullet)=K_0/R=\St_G$, the vanishing in positive degrees
proves~\textup{(1)}. The displayed isomorphism is induced by the inclusion
$\Cc(\Omega)\hookrightarrow K_0$, namely extension by zero. Composing
with $\varepsilon$ gives exactly the map in~\textup{(2)}. It is
$B$-equivariant because $\Omega$ is $B$-stable. Finally, the stabilizer
of $Bw_0$ in $B$ is $M$, so $\Omega\simeq M\backslash B\simeq U$.
In the coordinates $u\mapsto Bw_0u$ the actions are
\[
 (\rho(u_0)f)(u)=f(uu_0),\qquad
 (\rho(m)f)(u)=f(m^{-1}um)\quad (u_0,u\in U,\ m\in M).
\]
\end{proof}

\begin{remark}
The same open-cell map with integral coefficients identifies
$\St_{G,\mathbb Z}|_B$ with $\Cc(\Omega,\mathbb Z)$, where the integral
Steinberg representation is defined by the same parabolic quotient.
In particular it is free over $\mathbb Z$ \cite[Lemma~2.2]{BS76}.
\end{remark}

\subsection{Twisted coinvariants of the Steinberg representation}
\begin{corollary}\label{cor:steinberg-conjecture}
Let $\mathbf G$ be a connected reductive group over a finite or nonarchimedean local field $F$, and put $G=\mathbf G(F)$. If $F$ is local of positive characteristic, assume that its characteristic is sufficiently large in terms of the absolute root datum of $\mathbf G$. Then, for every smooth connected unipotent $F$-subgroup $\mathbf N\subseteq\mathbf G$ and every smooth character $\psi:N=\mathbf N(F)\to\C^\times$,
\[
 (\St_G)_{N,\psi}\ne\{0\}.
\]
In particular, when $\mathbf G$ is quasi-split, Conjecture~\ref{conj:degenerate} holds for the Steinberg representation with $a=1$ and without the exceptional case involving a vanishing Jacquet module. There is no characteristic restriction over finite fields.
\end{corollary}
\begin{proof}
Under the stated hypotheses on $F$, $N$ is conjugate into $U$. Over perfect fields this is the usual unipotent-subgroup theorem. In sufficiently large positive characteristic, apply Gille's unipotent-subgroup theorem \cite{Gi02} to the image of $N$ in each almost simple factor of the adjoint group, over its field of definition. The resulting parabolic pulls back to $\mathbf G$, and the central projection identifies their unipotent radicals. After this conjugation, Frobenius reciprocity and Theorem~\ref{thm:steinberg} give
\[
 \operatorname{Hom}_N(\St_G,\psi)
 \simeq\operatorname{Hom}_N(\Cc(U),\psi)\ne\{0\}.
\]
Explicitly, a nonzero functional is $f\mapsto\int_N f(n)\psi(n)^{-1}\,dn$, with counting measure when $F$ is finite.
\end{proof}

\begin{remark}[Containment of connected unipotent subgroups]\label{rem:unipotent-containment}
Every smooth connected $F$-split unipotent subgroup is contained in a $G$-conjugate of $U$, in every characteristic. Over a perfect field every smooth connected unipotent group is $F$-split. Connectedness alone does not suffice over a local field of positive characteristic, and unconditional nonvanishing for all $(N,\psi)$ is false there.

For example, take $F=\mathbb F_2((t))$, $L=F(\sqrt t)$ and $G=\mathrm{PGL}_2(F)$. Multiplication on the two-dimensional $F$-space $L$ embeds
\[
 \mathbf N=(\operatorname{Res}_{L/F}\mathbf G_m)/\mathbf G_m
 \ \subset\ \mathrm{PGL}_2,\qquad
 N=\mathbf N(F)=L^\times/F^\times.
\]
This is the standard smooth connected wound unipotent subgroup; see \cite[Example~4.19]{BMRS24}. Indeed, $L\otimes_F\overline F\simeq\overline F[\epsilon]/(\epsilon^2)$, so $\mathbf N_{\overline F}\simeq\mathbf G_a$. Explicitly,
\[
 N=\Set{\begin{bmatrix}a&tb\\b&a\end{bmatrix}:(a,b)\in F^2\setminus\{(0,0)\}},
\]
where brackets denote classes in $\mathrm{PGL}_2(F)$. Since $t$ is not a square in $F$, all displayed matrices are invertible. The action on $\mathbb P^1(F)$ is simply transitive: the image of $[1:0]$ is $[a:b]$. Thus $N$ is compact and has no fixed point in $\mathbb P^1(F)$, whereas every conjugate of $U$ fixes a point. In particular, $N$ is not conjugate into $U$.

For this rank-one group,
\[
 \St_G|_N\simeq C^\infty(N)/\C\mathbf1,
 \qquad (\St_G)_{N,1}=0.
\]
The second assertion follows from the first. This example shows that Corollary~\ref{cor:steinberg-conjecture} fails if the unipotent group is not $F$-split.
\end{remark}

\section{Grothendieck groups for restrictions to \texorpdfstring{$B$, $B_0$, and $U$}{B, B0, and U}}\label{sec:grothendieck-B}

Let $\mathbf G$ be a connected reductive group over $F$. Fix a maximal
$F$-split torus $\mathbf A$ in $\mathbf G$, let
$\mathbf M=Z_{\mathbf G}(\mathbf A)$ be the corresponding minimal Levi
subgroup, and fix a minimal parabolic subgroup
$\mathbf B=\mathbf M\mathbf U$. Write $G,A,M,B,U$ for the respective
groups of $F$-points, and put $B_0=AU$.

In the next few sections we discuss the restriction of an irreducible
representation of $G$ to $B$, $B_0$, or $U$. The resulting formulas
involve alternating sums, for which we need an appropriate Grothendieck group.
Let $H\in\{B,B_0,U\}$. Note that the Grothendieck group of all smooth
$H$-modules is zero.
Indeed, for any such module $V$, the direct sum
$W=\bigoplus_{j\geq 1}V$ satisfies $W\simeq V\oplus W$, whence
$[V]=0$. On the other hand, the restriction of a finite-length
admissible representation of $G$ to $B_0$ is finitely generated
(Lemma~\ref{lem:AU-finiteness}), and hence so is its restriction to $B$.
Its twisted $U$-coinvariants are also finite-dimensional for every
smooth character of $U$ (Lemma~\ref{lem:AU-finiteness} below).

Even for finitely generated modules there is a further issue. Consider
the exact sequence of modules over $R=\C[x,x^{-1}]$,
\[
 0\longrightarrow R\xrightarrow{\,x-1\,}R
 \xrightarrow{\,f\mapsto f(1)\,}\C\longrightarrow 0.
\]
The first map is multiplication by $x-1$. Thus the class of $\C$ in the
Grothendieck group of finitely generated $R$-modules is zero.

The same phenomenon occurs over a nonarchimedean local field when the
center of $\mathbf G$ contains a positive-dimensional $F$-split torus.
Put $Z=Z(\mathbf G)(F)$ and $Z_A=Z\cap A$. For $H=B$ or $B_0$,
if $V$ is a finitely generated smooth $H$-module on which $Z_A$
acts by a character, then
$[V]=0$ in the Grothendieck group of all finitely generated smooth
$H$-modules with unrestricted central character. Indeed, the valuation of a suitable $F$-rational
character of $\mathbf G$ gives a homomorphism $\nu:H\to\mathbb Z$
nonzero on $Z_A$. On $V\otimes R$ let $h$ act by
$h(v\otimes f)=hv\otimes x^{\nu(h)}f$. This module is finitely
generated: if $z\in Z_A$ has $\nu(z)=d>0$, translates by powers
of $z$ shift the exponent by multiples of $d$ and act on $V$ by
scalars. Thus generators of $V$ tensored with
$1,x,\ldots,x^{d-1}$ generate $V\otimes R$. Multiplication by $x-1$
now gives an exact sequence
\[
 0\longrightarrow V\otimes R\xrightarrow{\,x-1\,}V\otimes R
 \longrightarrow V\longrightarrow 0,
\]
which gives the asserted vanishing.

Keeping the above issues in mind, we will use categories of modules
for $B$, $B_0$, and $U$ that arise through restriction from finite-length
representations of $G$. This gives the finiteness properties we need,
though it has the drawback that the categories depend on the ambient
group $G$.
For $H\in\{B,B_0,U\}$, let $\mathcal R_G(H)$ be the full subcategory
of smooth $H$-modules consisting of the subquotients of $\pi|_H$,
where $\pi$ ranges over finite-length smooth admissible
representations of $G$. Thus its objects are the modules isomorphic
to $E_1/E_2$, with $E_2\subseteq E_1\subseteq\pi|_H$ both $H$-submodules.

These categories are closed under subobjects, quotients and finite direct
sums, the latter because we can take the direct sum of the ambient
$G$-representations. They are abelian, with short exact sequences
inherited from smooth $H$-modules. We do not take the closure under
extensions in the ambient category or under countable direct sums.
For a smooth character $\omega$ of $Z$, let
$\mathcal R_G^\omega(H)$, for $H=B$ or $B_0$, consist of the objects
on which $Z\cap H$ acts by $\omega|_{Z\cap H}$. We denote the
Grothendieck groups of $\mathcal R_G^\omega(B)$,
$\mathcal R_G^\omega(B_0)$ and $\mathcal R_G(U)$ respectively by
$K_0^\omega(B)$, $K_0^\omega(B_0)$ and $K_0(U)$, omitting from this
notation the fixed ambient group $G=\mathbf G(F)$. The prescribed
central action is on $Z$ for $B$ and on $Z_A$ for $B_0$.
Restriction preserves subquotients of the same $G$-representation
and is exact, giving maps
\[
  K_0^\omega(B)\longrightarrow K_0^\omega(B_0)\longrightarrow K_0(U).
\]
Over a finite field these are the usual categories of finite-dimensional
complex representations with the indicated central action: the
regular representation of $G$ restricts to a direct sum of copies
of the regular representation of $H$.

For each object of $\mathcal R_G(H)$, its twisted $U$-coinvariants are
finite-dimensional. This follows from the same property for the ambient
restriction (Lemma~\ref{lem:AU-finiteness}) and the exactness of
twisted coinvariants \cite[Proposition~1.9(a)]{BZ77}.

\begin{lemma}[Membership of the inducing modules]\label{lem:inducing-membership}
If $V$ is a closed subgroup of $U$, then $\cind_V^U\one$ belongs to
$\mathcal R_G(U)$. If $A$ normalizes $V$ and $\chi$ is a smooth
character of $A$, inflated to $AV$, then $\cind_{AV}^{B_0}\chi$ belongs
to $\mathcal R_G(B_0)$, with central action $\chi|_{Z_A}$. When
$B=TU$ with $T$ a torus, the analogous statement holds with $T$
and $B$ in place of $A$ and $B_0$.
\end{lemma}

\begin{proof}
By Theorem~\ref{thm:steinberg}, $\St_G|_U\simeq\Cc(U)$.
Integration along the left $V$-cosets gives the surjection
\[
 \Cc(U)\longrightarrow \cind_V^U\one,\qquad
 f\longmapsto\left(u\longmapsto\int_V f(vu)\,dv\right).
\]

Next, for any smooth character $\eta$ of $A$,
$\cind_A^{B_0}\eta$ is a subquotient of a principal series restriction.
To see this, choose a compact open subgroup $K_M$ of $M$ with
$A\cap K_M\subseteq\ker\eta$. Since $A$ is central in $M$, extend
$\eta$ trivially across $K_M$ to $AK_M$ and induce to $M$. The induced
module is finite-dimensional, since $M/A$ is compact, and $A$ acts
by $\eta$. Transport this coefficient by the longest Weyl element
so that it is the coefficient on the open Bruhat cell. The open-cell
module, restricted to $B_0$, is a finite direct sum of copies of
$\cind_A^{B_0}\eta$ (see Section~\ref{sec:principal}).

Write $m_V(a)$ for the scaling of Haar measure under conjugation by $a$,
so that $d(ava^{-1})=m_V(a)\,dv$. The same integration map, now on
functions on $B_0$, gives a surjection
\[
  \cind_A^{B_0}(\chi m_V^{-1})\longrightarrow\cind_{AV}^{B_0}\chi .
\]
The assertion for $B=TU$ follows in the same way.
\end{proof}

\Needspace{5\baselineskip}
\begin{lemma}[Twisted coinvariants of compact induction]\label{lem:compact-twisted-coinvariants}
Let $U$ be a unimodular totally disconnected locally compact group and
$J$ a closed unimodular subgroup. Let $E$ be a smooth $J$-module and
$\psi$ a smooth character of $U$. A choice of nonzero invariant measure on
$J\backslash U$ gives an isomorphism
\begin{equation}\label{eq:compact-twisted-coinvariants}
  (\cind_J^U E)_{U,\psi}\simeq E_{J,\psi|_J}.
\end{equation}
\end{lemma}

\begin{proof}
Write $[e]$ for the image of $e\in E$ in $E_{J,\psi|_J}$. For
$f\in\cind_J^U E$, the function $u\mapsto\psi(u)^{-1}[f(u)]$ descends
to a compactly supported function on $J\backslash U$. Integration
against the chosen measure therefore gives the map in
\eqref{eq:compact-twisted-coinvariants}. Its algebraic dual is the
Frobenius reciprocity isomorphism
\[
 \operatorname{Hom}_J(E,\C_{\psi|_J})
 \simeq\operatorname{Hom}_U(\cind_J^U E,\C_\psi)
\]
obtained by smooth duality and ordinary induction
\cite[Proposition~1.9(b),(d)]{BZ77}. Hence the map is an isomorphism.
\end{proof}

In our applications $U$ is unipotent. Every element of $U$ lies in a
compact subgroup, so $U$ and each of its closed subgroups are unimodular.

For later use, let $P=LV$ be a parabolic subgroup of $G$.
The contribution of a Bruhat cell to the restriction of
$\Ind_P^G\sigma$ to $H\in\{B,B_0,U\}$ is obtained by inflating the
coefficient through the Levi projection and then compactly inducing
to $H$. Both operations are exact. Thus if the coefficient is a
subquotient of the restriction of a finite-length admissible
representation $\sigma$ of the Levi, the resulting module is a
subquotient of the corresponding Bruhat term of
$(\Ind_P^G\sigma)|_H$. Since parabolic induction
preserves finite length and admissibility, this puts the image in
$\mathcal R_G(H)$. When a central character is prescribed, the
central action on the cell is read off from the Levi projection.
We use this observation to transport Grothendieck-group relations
from a Levi subgroup.

\section{Restriction to a minimal parabolic: principal series}\label{sec:principal}

Use the general connected reductive group and the minimal parabolic $B=MU$ of Section~\ref{sec:steinberg}. 
Let $\sigma$ be a finite-dimensional smooth complex representation of $M$ with a central character, inflated to $B$. Write $\omega$ for the restriction of its central character to $Z(\mathbf G)(F)$. Let $\delta=\delta_B$ be the modulus character of $B$ when $F$ is local, and the trivial character of $B$ when $F$ is finite. Define the un-normalised principal series
\[
\Ind_B^G( \sigma)
=\{f:G\to\sigma\text{ smooth}:f(bg)=\sigma(b)f(g), \, b\in B, \, g\in G\},
\]
on which $G$ acts by right translation.
For $w\in W$, put
\[
 B_w:=B\cap w^{-1}Bw,\qquad
 \sigma_w(h):=\sigma(whw^{-1}), \quad h\in B_w.
\]
The expression for $\sigma_w$ uses its inflation to $B$. 


\begin{proposition}\label{prop:principal-restriction}
The restriction $I(\sigma)|_B$ has a finite support filtration, indexed by $W$ in reverse Bruhat order, whose quotients are
\begin{equation}\label{eq:restriction-4-1}
 \gr_w(\Ind^G_B(\sigma)|_B)\simeq
 \cind_{B_w}^B\bigl(\sigma_w \bigr).
\end{equation}
In particular, in $K_0^{\omega}(B)$ we have
\begin{equation}\label{eq:restriction-4-2}
 [\Ind^G_B(\sigma)|_B]=\sum_{w\in W}
 [\cind_{B_w}^B(\sigma_w)].
\end{equation}
\end{proposition}

\begin{proof}
Filter smooth sections on $B\backslash G$ by compact support in the upward-closed open unions of cells. These are sections of the locally trivial vector bundle associated with 
$\sigma$. Restriction to a closed stratum gives the corresponding compactly supported sections: surjectivity follows by a finite compact open cover of the support on which the bundle is trivial and extending the coefficient functions in each trivialization by the compact-open extension lemma \cite[\S5.11]{BZ77}. The strata are $B_w\backslash B$. Explicitly, putting $\varphi(b)=f(wb)$ gives
\[
 \varphi(hb)=f((whw^{-1})wb)
 = \sigma_w(h)\varphi(b),\qquad h\in B_w.
\]
This identifies the cell module with the compact induction in
\eqref{eq:restriction-4-1}. Since $B=B_wB_0=B_wU$, restriction of the
cell module to $B_0$ or $U$ gives compact induction from $B_w\cap B_0$
or $B_w\cap U$, respectively. All the filtration terms and their
successive quotients belong to $\mathcal R_G(B)$, since they are
subquotients of the principal series restricted to $B$. The center
acts by $\omega$ throughout. Adding the short exact sequence
relations gives~\eqref{eq:restriction-4-2}.
\end{proof}

\section{Restriction to the Borel subgroup of \texorpdfstring{$\GL_n(F)$}{GL(n,F)}} \label{sec:GL}

Let $G_n=\GL_n(F)$ and $Z_n\simeq F^\times$ be its center. Let $B_n=T_n U_n$ be the  upper triangular Borel subgroup of $G_n$, where $T_n$ is the diagonal torus and $U_n$ is as before the upper triangular maximal unipotent subgroup of $G_n$. 
In this section, we will prove the following result.

\begin{theorem}\label{thm:GL}
For every irreducible smooth representation $\pi$ of $G_n$ with central character $\omega$, we have 
\begin{equation} \label{GLsum}
[\pi|_{B_n}]= \sum^r_{i=1} n_i [\cind^{B_n}_{B_{w_i}}\chi_i]\in K_0^{\omega}(B_n),
\end{equation}
where $n_i\in\mathbb Z$, $w_i\in W\simeq S_n$ is a Weyl group element, and $\chi_i: T_n\to \mathbb{C}^\times$ is a character extending $\omega$ on $Z_n$, inflated to $B_{w_i}$, $i=1,2,\ldots, r$.
\end{theorem}

\subsection{Cuspidal representations}
We first prove Theorem \ref{thm:GL} when $\pi$ is cuspidal. Write
\[
P_n= G_{n-1}\ltimes V_n= \left\{\begin{pmatrix}g&v\\0&1\end{pmatrix}:
 g\in G_{n-1},\ v\in F^{n-1}\right\}
\]
for the mirabolic subgroup introduced in Section~\ref{sec:extension}, where $V_n = F^{n-1}$ is the unipotent radical of $P_n$. Recall the fixed nontrivial additive character $\psi_F$ of $F$. Define a character
\[
\psi_n: U_n\to\mathbb C^\times, \qquad (u_{ij})\mapsto \psi_F(u_{12}+u_{23}+\cdots+u_{n-1,n}).
\]
For every irreducible cuspidal representation $\pi$ of $G_n$, the mirabolic restriction theorem for nonarchimedean local fields \cite[\S5.18]{BZ76} and its finite-field counterpart \cite[Theorem~2.3 and \S3]{Gel70} give
\begin{equation}\label{eq:restriction-5-1}
\pi|_{P_n}\simeq \cind_{U_n}^{P_n}\psi_n.
\end{equation}

Let $L_n = B_{n-1} Z_n\subset G_n$, where $B_{n-1}$ is viewed as a subgroup of $G_n$ via the embedding 
\[
g\longmapsto \begin{bmatrix} g \\ & 1\end{bmatrix}.
\]
Then $B_n = L_n \ltimes V_n$. 
Define a representation $\tau_n$ of the parabolic subgroup 
\[
Q_n:=P_nZ_n \subset G_n
\]
by
\[
\tau_n: = \left(\cind_{U_n}^{P_n}\psi_n\right)\boxtimes \mathbf1_{Z_n}.
\]
Let $\delta_{Q_n}$ be the modulus character of $Q_n$ when $F$ is local, and the trivial character of $Q_n$ when $F$ is finite.

\begin{lemma} \label{lem:cusp}
    For $n\ge2$, there is a short exact sequence of $B_n$-modules
    \[
    0 \longrightarrow \tau_n|_{B_n} \longrightarrow \cind^{B_n}_{L_n}(\tau_{n-1}|_{B_{n-1}}\boxtimes \mathbf1_{Z_n}) \longrightarrow {\rm Inf}^{B_n}_{L_n}\left( (\tau_{n-1}|_{B_{n-1}}\boxtimes \mathbf1_{Z_n}) \otimes \delta_{Q_n}|_{L_n}\right)\longrightarrow 0,
    \]
    where ${\rm Inf}$ denotes inflation of a representation of $L_n$ to $B_n$.
\end{lemma}

\begin{proof}
    We first work on 
$P_n = G_{n-1}\ltimes V_n =  G_{n-1}\ltimes F^{n-1}$.
    Choose an auxiliary irreducible cuspidal representation $\sigma$ of
$G_{n-1}$ with trivial central character, so that its restriction to
$B_{n-1}$ is precisely $\tau_{n-1}|_{B_{n-1}}$.\footnote{For $d=n-1>1$
and $k=\mathbb F_q$, a faithful character of
$\mathbb F_{q^d}^\times/k^\times$ has Frobenius orbit of length $d$
since its order $(q^d-1)/(q-1)$ exceeds $q^{d-1}-1$. This proves that there are always
cuspidal representations of $\GL_d(k)$ with trivial central character. Over a local
field with residue field $k$, inflate this representation to
$\GL_d(\mathcal O_F)$, extend trivially across scalar uniformizers
and compactly induce from $F^\times\GL_d(\mathcal O_F)$, the usual
depth-zero construction. For $d=1$ take the trivial character.}  We have the realization 
    \[
    \cind^{P_n}_{G_{n-1}}\sigma = C^\infty_c(F^{n-1}, \sigma),
    \]
    with the group action given by
    \[
    (g\cdot f)(x) = \sigma(g)f(g^{-1}x),\qquad (v\cdot f)(x) = f(x+v),
    \]
where $f\in C^\infty_c(F^{n-1},\sigma)$, $g\in G_{n-1}$ and $v, x\in F^{n-1}$. The integration over $F^{n-1}$ defines a $P_n$-equivariant surjection 
\[
\phi: \cind^{P_n}_{G_{n-1}}\sigma  \to {\rm Inf}^{P_n}_{G_{n-1}}(\sigma|\det|_F),\qquad f\mapsto \int_{F^{n-1}}f(x)dx.
\]
Here and henceforth, $|\cdot|_F$ is the normalized absolute value when $F$ is local, and is by convention the trivial character of $F^\times$ when $F$ is finite. 

Recall the dual space $F_{n-1}$ from Section~\ref{sec:cyclic}. Fix a Haar
measure $dx$ on $F^{n-1}$ (counting measure when $F$ is finite), and put
\[
 \widehat f(\lambda)=\int_{F^{n-1}}f(x)\psi_F(-\lambda x)\,dx,
 \qquad \lambda\in F_{n-1}.
\]
Then
\[
 \widehat{g\cdot f}(\lambda)
   =|\det g|_F\sigma(g)\widehat f(\lambda g),\qquad
 \widehat{v\cdot f}(\lambda)
   =\psi_F(\lambda v)\widehat f(\lambda).
\]
Moreover, $\phi(f)=\widehat f(0)$, so Fourier transform identifies
$\ker\phi$ with $C_c^\infty(F_{n-1}\setminus\{0\},\sigma)$.
The stabilizer of $\lambda_0=(0,\ldots,0,1)$ in $G_{n-1}$ is
$P_{n-1}$, and its action on the fiber at $\lambda_0$ is
$\sigma|_{P_{n-1}}\otimes|\det|_F$. Consequently,
\[
 \ker\phi\simeq
 \cind_{P_{n-1}\ltimes V_n}^{P_n}
 \left((\sigma|_{P_{n-1}}\otimes|\det|_F)
              \otimes\psi_n|_{V_n}\right).
\]
The determinant twist can be removed from this coefficient:
\[
 \sigma|_{P_{n-1}}\otimes|\det|_F
 \simeq
 \left(\cind_{U_{n-1}}^{P_{n-1}}\psi_{n-1}\right)
                    \otimes|\det|_F
 \simeq \cind_{U_{n-1}}^{P_{n-1}}\psi_{n-1},
\]
because $|\det|_F$ is trivial on $U_{n-1}$.
Here the first isomorphism is \eqref{eq:restriction-5-1} in rank $n-1$;
for $n=2$ both lower-rank groups are trivial. Induction in stages now gives
\[
 \ker\phi\simeq
 \cind_{P_{n-1}\ltimes V_n}^{P_n}
       \left(\cind_{U_n}^{P_{n-1}\ltimes V_n}\psi_n\right)
 \simeq\cind_{U_n}^{P_n}\psi_n.
\]
The lemma follows from taking the tensor product of the short exact sequence 
\[
0\longrightarrow \cind^{P_n}_{U_n}\psi_n \longrightarrow \cind^{P_n}_{G_{n-1}}\sigma  \longrightarrow {\rm Inf}^{P_n}_{G_{n-1}}(\sigma|\det|_F) \longrightarrow 0
\]
with $\mathbf1_{Z_n}$ and the restriction to $B_n$, and using \eqref{eq:restriction-5-1} again.
\end{proof}

For $J=\{j_1<\cdots<j_s\}\subseteq\{2,\ldots,n\}$, write
$\{2,\ldots,n\}\setminus J=\{k_1<\cdots<k_r\}$ and define
$w_J\in S_n$ by
\[
 \bigl(w_J^{-1}(1),\ldots,w_J^{-1}(n)\bigr)
       =(k_r,\ldots,k_1,1,j_1,\ldots,j_s).
\]
For $i<j$, this ordering gives
\[
 w_J(i)<w_J(j)\quad\Longleftrightarrow\quad j\in J.
\]
It follows, by checking the matrix root subgroup
$\{1_n+xE_{ij}:x\in F\}$ for each $i<j$, that
\[
 B_{w_J}=B_n\cap w_J^{-1}B_nw_J=T_nV_J,\qquad
 V_J=\{(u_{ij})\in U_n:u_{ij}=0\text{ if }i<j\text{ and }j\notin J\}.
\]
In particular, $w_\varnothing=w_0$ and $w_{\{2,\ldots,n\}}=1$.
Define a character 
\[
 \delta_J : T_n \to\mathbb C^\times, \qquad  {\rm diag}(t_1,\ldots,t_n)\mapsto \prod_{j\in J}\prod_{i<j}|t_i/t_j|_F
\]
when $F$ is local. By convention, set $\delta_J$ to be the trivial character when $F$ is finite. 
For a character $\omega$ of $F^\times$, define 
\begin{equation} \label{lambdaomega}
\lambda_\omega: T_n\to \mathbb C^\times,\qquad {\rm diag}(t_1,\ldots, t_n)\mapsto \omega(t_n).
\end{equation}

We have the following special case of Theorem \ref{thm:GL}.

\begin{proposition}\label{prop:GL}
    Let $\pi$ be an irreducible cuspidal representation of $G_n$ with central character $\omega$. Then 
    \[
    [\pi|_{B_n}] =[\tau_n|_{B_n}\otimes \lambda_\omega]  = \sum_{J\subset\{2,\ldots, n\}}(-1)^{|J|}\cind^{B_n}_{B_{w_J}}(\lambda_\omega \delta_J) \in  K_0^{\omega}(B_n).
    \]
\end{proposition}

\begin{proof}
    By \eqref{eq:restriction-5-1}, we have 
    \[
    \pi|_{B_n} \simeq \tau_n|_{B_n}\otimes \lambda_\omega.
    \]
We check that the exact sequence of Lemma~\ref{lem:cusp} and its
lower-rank substitutions can be used in the new categories. For a
$B_{n-1}$-module $E$, let $E^{(\omega)}$ be its extension to
$L_n=B_{n-1}Z_n$ on which $Z_n$ acts by $\omega$. Suppose that
$E$ is a subquotient of $\rho|_{B_{n-1}}$ for a finite-length
admissible representation $\rho$ of $G_{n-1}$. On the Levi
$G_{n-1}\times G_1$ of $Q_n$, use the representation
\[
       (g,a)\longmapsto \omega(a)\rho(a^{-1}g).
\]
It is finite-length admissible, with $Z_n$ acting by $\omega$.
The open cell for induction from the opposite parabolic to $Q_n$,
with this Levi, has stabilizer $L_n$ in $B_n$. Its coefficient is
the restriction of the displayed Levi representation to $L_n$.
By exactness, $\cind_{L_n}^{B_n}E^{(\omega)}$ is a subquotient of
this open-cell module. In particular its ambient $G_n$-representation
has central character $\omega$.

In the realization on $\Cc(V_n,E^{(\omega)})$, integration over $V_n$
gives the surjection onto
\[
 \operatorname{Inf}_{L_n}^{B_n}
       \bigl(E^{(\omega)}\otimes\delta_{Q_n}|_{L_n}\bigr).
\]
Indeed, the action of $L_n$ on the integral is twisted by the
scaling of Haar measure under its action on $V_n$, namely
$\delta_{Q_n}|_{L_n}$. Thus the inflated term also belongs to the
category, with the same ambient representation. Both operations
on $E$ are exact. Apply them to the lower-rank relations in the
sequence of Lemma~\ref{lem:cusp}, taking the $Z_n$-action to be
$\omega$. Its kernel $\tau_n|_{B_n}\otimes\lambda_\omega$ is also
a submodule of the middle term, with the same ambient representation.
Iterating the sequence from $\tau_1=\one$ gives the displayed
alternating sum, with the twists $\lambda_\omega\delta_J$ supplied
by the integration factors.
\end{proof}

\subsection{Bernstein-Zelevinsky filtration} We now prove Theorem~\ref{thm:GL} for arbitrary irreducible smooth representations $\pi$ of $G_n$ by an inductive argument.

Let $\omega$ be the central character of $\pi$. The local mirabolic filtration is \cite[\S3.5]{BZ77}. Over a finite field the corresponding decomposition follows by iterating the two orbits of $G_{r-1}$ on the character group of $V_r$: the zero character gives a $G_{r-1}$-module, and a nonzero character has stabilizer $P_{r-1}$; see \cite[\S1]{Gel70}. Thus $\pi|_{Q_n}= \pi|_{P_nZ_n}$  has a filtration with subquotients of the form 
\[
I_k(\sigma):=\cind^{P_n}_{P_{n,k}}(\sigma\boxtimes \psi_{n-k}) \boxtimes \omega,
\]
where $k=0,1,\ldots, n-1$, $\sigma\in {\rm Irr}(G_k)$ and 
\[
P_{n,k}:=\Set{\begin{bmatrix} g & X \\ & u\end{bmatrix}: g\in G_k, \ u\in U_{n-k}, \ X\in F^{k\times (n-k)}}.
\]
Here the representation $\sigma\boxtimes \psi_{n-k}$ of $G_k\times U_{n-k}$ is inflated to $P_{n,k}$.

If $k=0$, then by Proposition \ref{prop:GL} we have 
\[
[I_0(\mathbf 1)|_{B_n}] = [\tau_n|_{B_n}\otimes \lambda_\omega] \in  K_0^{\omega}(B_n)
\]
is of the form \eqref{GLsum}.

Next assume that $1\leq k<n$. Let $Q$ be the standard parabolic subgroup of type $(k,n-k-1,1)$. We identify the symmetric group $S_n$ with the group of permutation matrices in 
$G_n$. For every minimal length representative $w\in S_{n-1}\subset S_n$ of 
\[
(S_k \times S_{n-k-1})\backslash S_{n-1},
\]
put 
\[
H_w:= B_n \cap w^{-1}Qw. 
\]
The composition 
\[
H_w \xrightarrow{{\rm Ad}_w} wB_nw^{-1} \cap Q \hookrightarrow Q \twoheadrightarrow G_k\times Q_{n-k}
\]
maps onto $B_k\times B_{n-k}$, which gives a quotient homomorphism
\begin{equation} \label{quo}
H_w \longrightarrow B_k\times B_{n-k}.
\end{equation}
Then $I_k(\sigma)$ has a filtration with subquotients of the form
\begin{equation} \label{Hw}
\cind^{B_n}_{H_w} \left(\sigma|_{B_k} \boxtimes \left(\tau_{n-k}|_{B_{n-k}}\otimes \lambda_{\omega \cdot \omega_\sigma^{-1}}\right)\right),
\end{equation}
where $\lambda_{\omega\cdot\omega_\sigma^{-1}}: T_{n-k}\to \mathbb{C}^\times$ is given by \eqref{lambdaomega} and $\omega_\sigma$ denotes the central character of $\sigma$, and the inducing representation is inflated 
to $H_w$ through \eqref{quo}.
Apply the induction hypothesis to the first factor and
Proposition~\ref{prop:GL} to the second, keeping both block-central
characters fixed. Since $w$ fixes the last index, $H_w$ is also the
intersection of $B_n$ with the conjugate of the two-block parabolic of
type $(k,n-k)$, and \eqref{quo} is its Levi projection. If the two
coefficients are subquotients of restrictions of finite-length
admissible representations of $G_k$ and $G_{n-k}$, their external
tensor product is a subquotient of the restriction of the corresponding
external tensor product on the Levi. The Bruhat-cell observation
of Section~\ref{sec:grothendieck-B} therefore shows that inflation
and compact induction in~\eqref{Hw} take these coefficients into
$\mathcal R_{G_n}(B_n)$. The operations are exact, so the two block
relations and the Bruhat filtration give relations in the prescribed
central-character category. The ambient representations used in
these relations can be taken with the stated central characters:
for the cuspidal terms this was shown in the proof of
Proposition~\ref{prop:GL}, and in the inductive step the two block
characters multiply to $\omega$ on $Z_n$.
By compact induction in stages, the resulting terms are
induced from the inverse images in $H_w$ of the two block stabilizers.
If $v\in S_k\times S_{n-k}$ indexes the two block stabilizers,
their inverse image is $B_{vw}=B_n\cap(vw)^{-1}B_n(vw)$.
Every term has character $\omega$ on $Z_n$, which proves
\eqref{GLsum} in $K_0^\omega(B_n)$.

This finishes the proof of Theorem \ref{thm:GL}.

\section{Restriction conjectures}\label{sec:question}

The aim of this section is to formulate certain conjectures on representations
of $G$ restricted to $B_0$ and to $U$. Before doing so, we specify the class of representations
to which the conjectures apply. For finite reductive groups, these
conjectures describe certain relationships among character values at unipotent
elements; see \cite{NP} for related work. Over a nonarchimedean local
field there is no comparable pointwise character theory on the
unipotent set, so there is no similarly direct interpretation.

\subsection{Stable representations}

For a nonarchimedean local field $F$, let $K_0(G)$ denote the Grothendieck group of
finite-length smooth admissible representations of $G$, and put
$K_0(G)_\C:=K_0(G)\otimes_{\mathbb Z}\C$. An element
$\Pi=\sum_i c_i[\pi_i]\in K_0(G)_\C$, with finitely many $c_i\in\C$ and $\pi_i\in{\rm Irr}(G)$, is called
\emph{stable} if its Harish--Chandra character
$\Theta_\Pi=\sum_i c_i\Theta_{\pi_i}$ is constant on stable conjugacy
classes of strongly regular semisimple elements. 

When $F=\mathbb F_q$ is finite, let $K_{\mathrm{unif}}(G)$ be the
subgroup of the Grothendieck group of finite-dimensional complex
representations of $G$ generated by Deligne--Lusztig virtual representations. Put
\[
K_{\mathrm{unif}}(G)_\C:=K_{\mathrm{unif}}(G)\otimes_{\mathbb Z}\C.
\]
Thus
\[
K_{\mathrm{unif}}(G)_\C=
 \operatorname{Span}_{\C}\{R_{\mathbf T'}^{\mathbf G}(\theta):
 \mathbf T'\text{ a rational maximal torus},\
 \theta\in\operatorname{Hom}(\mathbf T'(F),\C^\times)\}
\]
is the space of \emph{uniform representations} of $G$ \cite{DL76}. Thus a uniform representation $\Pi$ is a complex linear combination of Deligne--Lusztig virtual representations. In good characteristic its character $\Theta_\Pi$ is \emph{weakly stable} in the following sense. If $\mathrm{Fr}$ is Frobenius and $g\in G$, choose $x\in\mathbf G(\overline{\mathbb F}_q)$
with $g=x^{-1}\mathrm{Fr}(x)$ and define the Shintani twisting
$\operatorname{Sh}(g)=\mathrm{Fr}(x)x^{-1}$. This is again an element of $G$ geometrically conjugate to $g$.
Then, by Digne--Michel \cite[IV, Corollary~2.5]{DM85},
\[
\Theta_\Pi\circ\operatorname{Sh}=\Theta_\Pi.
\]


\subsection{Grothendieck groups}
Fix a maximal $F$-split torus $\mathbf A$ in $\mathbf G$, put
$\mathbf M=Z_{\mathbf G}(\mathbf A)$, and let $\mathbf B=\mathbf M\mathbf U$
be a minimal parabolic subgroup containing $\mathbf M$. Write
$M:=\mathbf M(F)$, $B:=\mathbf B(F)$, and $U:=\mathbf U(F)$. Put
\[
 \mathbf B_0= \mathbf A \mathbf U, \qquad A:=\mathbf A(F), \qquad B_0:=\mathbf B_0(F)= AU \subset B.
\]
Then $\mathbf B_0$ is a maximal $F$-split connected solvable subgroup of
$\mathbf G$. Inside this fixed $\mathbf B$ it is the largest such
subgroup: its unipotent radical is $\mathbf U$, and the image of a
connected split solvable subgroup in $\mathbf M=\mathbf B/\mathbf U$
lies in its central split torus $\mathbf A$. We have
$\mathbf B_0=\mathbf B$ if and only if $\mathbf G$ is split.

We use $\mathcal R_G^\omega(B_0)$ and $\mathcal R_G(U)$ of
Section~\ref{sec:grothendieck-B}, and write
\[
 K_0^\omega(B_0)_\C:=K_0^\omega(B_0)\otimes_{\mathbb Z}\C,\qquad
 K_0(U)_\C:=K_0(U)\otimes_{\mathbb Z}\C.
\]
The next lemma records two finiteness properties of restrictions.

\begin{lemma}\label{lem:AU-finiteness}
If $\pi$ is a finite-length smooth admissible representation of $G$,
then $\pi|_{B_0}$ is finitely generated. Moreover,
$\dim\pi_{U,\psi}<\infty$ for every smooth character $\psi: U\to \C^\times$.
\end{lemma}

\begin{proof}
The quotient $M/A$ is compact, and Iwasawa decomposition gives
$G=B_0C$ for a compact subset $C$. Choose finitely many $G$-generators
of $\pi$. Smoothness implies that the $C$-orbit of each of them spans a
finite-dimensional space; bases of these spaces generate $\pi|_{B_0}$.
The second assertion is the finite-multiplicity theorem
\cite[Theorem~1.0.4 and Appendix~B]{ABS22}, followed by exactness of
twisted coinvariants to pass from irreducible to finite-length modules.
\end{proof}

\subsection{The two conjectures}\label{subsec:two-conjectures}
In the identity for $B_0$, a virtual representation with several
central characters is treated one central character at a time:
write $\Pi=\sum_\omega\Pi_\omega$ by collecting its irreducible
terms with central character $\omega$. 

\begin{conjecture} \label{conj:restriction-A}
Let $\Pi\in K_0(G)_\C$ be stable when $F$ is local and $\Pi\in K_{\mathrm{unif}}(G)_\C$ when $F$ is finite.
For each central character $\omega$ occurring in $\Pi$, there are finitely many closed subgroups $H_i$ with $A\subseteq H_i\subseteq B_0$, smooth characters
$\chi_i:A\to\C^\times$ agreeing with $\omega$ on $Z_A$ and inflated along $H_i\to A$, and $c_i\in\C$ such that 
\begin{equation}\label{Bres}
 [\Pi_\omega|_{B_0}]=\sum_i c_i [\cind^{B_0}_{H_i} \chi_i] \in K_0^{\omega}(B_0)_\C.
\end{equation}
\end{conjecture}

\begin{conjecture} \label{conj:restriction-B}
Let $\Pi\in K_0(G)_\C$ be stable when $F$ is local, and $\Pi\in K_{\mathrm{unif}}(G)_\C$ when $F$ is finite.
Then there are finitely many closed subgroups $V_i \subset U$ normalized by $A$ and $c_i\in\C$ such that 
\begin{equation}\label{eq:conjecture-B}
 [\Pi|_U]=\sum_i c_i [\cind^U_{V_i}\mathbf 1]\in K_0(U)_\C. 
\end{equation}
\end{conjecture}

By definition, $\pi|_{B_0}\in\mathcal R_G^\omega(B_0)$ and
$\pi|_U\in\mathcal R_G(U)$ for a finite-length smooth admissible
representation $\pi$ of $G$ with central character $\omega$. The inducing terms on the
right-hand sides of \eqref{Bres} and \eqref{eq:conjecture-B}
belong to these categories by Lemma~\ref{lem:inducing-membership}.

Restriction to $U$ gives the homomorphism
$K_0^\omega(B_0)\to K_0(U)$ of Section~\ref{sec:grothendieck-B}. For each term in \eqref{Bres},
$A\subseteq H_i$ gives $B_0=H_iU$, and $\chi_i$ is trivial on $H_i\cap U$.
Restriction of functions consequently gives
\[
  (\cind_{H_i}^{B_0}\chi_i)|_U\simeq\cind_{H_i\cap U}^{U}\one .
\]
The subgroup $H_i\cap U$ is normalized by $A$. Applying the homomorphism
to \eqref{Bres} and summing over $\omega$ proves that
Conjecture~\ref{conj:restriction-A} implies Conjecture~\ref{conj:restriction-B}. 

\begin{remark}\label{rem:exceptional-sl2}
By Theorem \ref{thm:GL},  every irreducible smooth representation of $\GL_n(F)$  satisfies \eqref{Bres}.  On the other hand, not all irreducible representations of $\mathrm{SL}_2(\mathbb F_q)$ with $q>3$ odd satisfy \eqref{Bres}: the exceptional irreducible representations of dimensions $(q+1)/2$ and $(q-1)/2$ do not satisfy it. Indeed, the only subgroups of $B$ containing $T$ are $T$ and $B$. Characters induced from $T$ vanish on $U\setminus\{1\}$, and one-dimensional characters of $B$ are trivial on $U$. Therefore the character of the representation
${\rm Ind}^B_H\chi$
is constant on $U\setminus\{1\}$, whereas these exceptional characters of dimensions $(q+1)/2$ and $(q-1)/2$ distinguish between the square and nonsquare parameters in $U$. 

Over a local field the individual constituents of a normalized
$\mathrm{SL}_2(F)$ principal series induced from a nontrivial quadratic
character give an analogous obstruction to B: their Whittaker
multiplicities vary with the nontrivial character of $U$. Take $\operatorname{char} F\ne2$. Every $A$-normalized subgroup of the
additive group $U=F$ is closed under multiplication by squares and
hence by their differences, which exhaust $F$. Thus $V=0$ or $F$,
and each $\cind_V^U\one$ has constant multiplicity on nontrivial
characters.

Thus stability is a natural condition for the restriction conjectures, and stable linear combinations of representations in an $L$-packet are the intended examples. The question is whether these representations of $G$ have restrictions to $B_0$ satisfying \eqref{Bres}.
\end{remark}

\begin{remark}[Linear dependence in Conjecture~\ref{conj:restriction-B}]\label{rem:restriction-linear-dependence}
Over a finite field, the compactly induced representations on the
right-hand side of \eqref{eq:conjecture-B} are not necessarily linearly
independent over $\mathbb Q$, even if the subgroups are restricted to the Weyl intersections
$U\cap wUw^{-1}$. Indeed, let $G=\GL_3(\mathbb F_q)$, let $U$ be upper unitriangular,
and put $X=U_{12}$, $Y=U_{23}$, and $Z=U_{13}=Z(U)$, where $U_{ij}$ is
the corresponding root subgroup. The four subgroups $X,Y,XZ,YZ$ are
distinct Weyl intersections and are normalized by the diagonal torus $A$.
Writing $P_H=\cind_H^U\one$, we have in $K_0(U)$
\[
 [P_X]-[P_{XZ}]=[P_Y]-[P_{YZ}].
\]

This follows by counting fixed cosets. Off the
center, the characters of $P_X$ and $P_{XZ}$ both have value $q$ if
the $(2,3)$-entry is zero, and both are zero otherwise. The analogous
statement for $P_Y$ and $P_{YZ}$ uses the $(1,2)$-entry. On $Z$ the
characters of $P_X$ and $P_Y$ have value $q^2$ at $1$ and zero elsewhere,
whereas those of $P_{XZ}$ and $P_{YZ}$ are identically $q$.

Thus the coefficients in \eqref{eq:conjecture-B} need not be unique.
The same example works over nonarchimedean local fields, with a
different proof that we do not give here.

It is an interesting question to find a family of $A$-invariant subgroups
of $U$ whose induced representations form a $\mathbb Q$-basis of the
subspace of $K_0(U)\otimes_{\mathbb Z}\mathbb Q$ spanned by the
representations $\cind_V^U\one$ for all $A$-invariant subgroups $V$ of $U$.
(Note that $U$ has only finitely many closed $A$-invariant subgroups if
$F$ is finite, or if $\mathbf G$ is split and $\operatorname{char}F\ne 2$.)
\end{remark}

The following lemma verifies a necessary separation property for
Conjecture~\ref{conj:restriction-B} over finite fields in good characteristic.

\begin{lemma}[Separation of geometric unipotent classes]\label{lem:separation-geometric-unipotent}
Suppose that $F=\mathbb F_q$ has good characteristic for $\mathbf G$.
If $u,v\in U$ are not conjugate in $\mathbf G(\overline{\mathbb F}_q)$,
there is an $F$-rational subgroup $\mathbf V\subseteq\mathbf U$,
normalized by the maximal torus $\mathbf T=Z_{\mathbf G}(\mathbf A)$,
such that the character of $\Ind_{\mathbf V(F)}^U\one$ is positive
at one of $u,v$ and zero at the other. In particular, the characters on
the right-hand side of \eqref{eq:conjecture-B} separate geometrically
distinct unipotent elements of $U$.
\end{lemma}

\begin{proof}
Interchange $u$ and $v$ if necessary so that
$v\notin\overline{\mathcal O_u}$, where $\mathcal O_u$ is the geometric
conjugacy class of $u$. Choose a rational Dynkin cocharacter
$\lambda\colon\mathbf G_m\to\mathbf G$ associated to $u$. The subgroup
$\mathbf U(\lambda)_{\geq 2}$ is generated over $\overline{\mathbb F}_q$ by
the root subgroups of $\lambda$-weight at least two. It is normal in
the associated parabolic $\mathbf P(\lambda)$, contains $u$, and meets
$\mathcal O_u$ in a dense open subset; in particular,
$\mathbf U(\lambda)_{\geq 2}\subseteq\overline{\mathcal O_u}$
\cite[Lemma~4.7 and \S6.1]{BDT20}; see also
\cite[Theorem~3.22 and \S4]{T16}.

By rational Bruhat decomposition, $\mathbf B$ and $\mathbf P(\lambda)$
contain a common rational maximal torus $\mathbf S$. The torus $\mathbf S$
normalizes $\mathbf U(\lambda)_{\geq 2}$ by normality in $\mathbf P(\lambda)$.
Choose $a\in U$ such that $a\mathbf S a^{-1}=\mathbf T$, and put
\[
\mathbf V=\mathbf U\cap a\mathbf U(\lambda)_{\geq 2}a^{-1}.
\]
The subgroup $\mathbf V$ is rational and $\mathbf T$-normalized. It
contains $aua^{-1}$ and is contained in $\overline{\mathcal O_u}$. Hence
$u$ fixes a coset in $U/\mathbf V(F)$, whereas $v$ fixes none. The
induced character counts these fixed cosets.
\end{proof}

\Needspace{9\baselineskip}
\begin{remark}\label{rem:separation-weyl-intersections}
For $\GL_n$ over any finite field and for split $G_2$ in
characteristic $p>3$, the characters induced from $V=U\cap U^w$
are sufficient to separate geometrically distinct unipotent elements. For $\GL_n$, this
already follows from the proof of Conjecture~\ref{conj:restriction-B}
with Weyl intersections (Theorem~\ref{thm:GL} after restriction to $U$),
since irreducible characters separate conjugacy classes and geometric
and rational unipotent conjugacy coincide. Here is a direct argument.

Let $u,v\in U$ have different Jordan types. After interchanging them
if necessary, choose $j$ with
$\dim\ker(v-1)^j<\dim\ker(u-1)^j$. Put $F_i=\ker(u-1)^i$ for
$0\leq i\leq n$. The Bruhat decomposition puts this flag in coordinate
position by $U$-conjugation: each $B$-orbit on the corresponding
partial flag variety contains a coordinate flag, and the diagonal
torus fixes that flag, so its $B$-orbit is already a $U$-orbit.
Replace $u$ and the flag by these conjugates.

For a standard basis vector $e_r$, let $b(r)$ be the first step of
the flag containing $e_r$. The subgroup
\[
 V=\{x\in U:(x-1)F_i\subseteq F_{i-1}\text{ for }1\leq i\leq n\}
\]
allows the entry in position $(r,s)$ exactly when $r<s$ and
$b(r)<b(s)$. To see that this is a Weyl intersection, order the
basis by increasing $b(r)$ and by decreasing $r$ within each block.
An entry above the diagonal in the standard order is also above the
diagonal in the new order precisely when $b(r)<b(s)$. Thus $V$ is
the intersection of the upper unitriangular groups for the two
orders.

By construction $u\in V$. For every $x\in V$, iterating the
flag-lowering condition gives $F_j\subseteq\ker(x-1)^j$, so no
conjugate of $v$ belongs to $V$. The character of $\Ind_V^U\one$
counts fixed cosets in $U/V$: it is positive at $u$ and zero at $v$.

For split $G_2$, the sum of the
Weyl-intersection characters is the permutation character on
$G/B$ by the Bruhat decomposition. Its values on the
seven rational unipotent classes are pairwise distinct, as seen in the
split-torus row of the Green function table \cite{ChevieJL}.
\end{remark}

\section{Some evidences for the conjectures}\label{sec:restriction-evidence}
We first check that Conjectures~\ref{conj:restriction-A}
and~\ref{conj:restriction-B} hold for a parabolically induced
representation if they hold for the inducing representation.

\subsection{Ordinary parabolic induction}\label{subsec:ordinary-parabolic}

Before proceeding further, we note the following well-known lemma, cf. \cite[Propositions 3.4.8 and 3.4.9]{DM20}.

\begin{lemma}\label{lem:parabolic-intersection}
Fix a maximal split torus $A$ in $G$. Let $B = MU$ be a minimal
parabolic subgroup in $G$ with $M=Z_G(A)$, and let $P=LV$ be a parabolic subgroup of $G$ with $L$ containing $M$. Then,
$M\cdot (L\cap U)$ is a minimal parabolic subgroup of $L$, and we have the Levi decomposition of $P\cap B$, which need not be a parabolic subgroup,
\[ P \cap B = M\cdot (L \cap U) \cdot (V \cap U).
\]
\end{lemma}

\begin{proposition}\label{prop:induction-permanence}
Let $\sigma$ be an irreducible admissible representation of a Levi
subgroup $L$ of a parabolic subgroup $P=LV$. If $\sigma$ satisfies
the restriction identity of Conjecture~\ref{conj:restriction-A} or
of Conjecture~\ref{conj:restriction-B}, then
$\pi=\cind_P^G\sigma=\Ind_P^G\sigma$ satisfies the corresponding
identity. More generally, the same implication holds for an arbitrary
finite-length admissible $L$-module $\sigma$, with \eqref{Bres} taken
central character by central character in its Grothendieck class.
\end{proposition}

\begin{proof}
Since $P\backslash G$ is compact, the two inductions agree.
Over a nonarchimedean local field, parabolic induction preserves
admissibility and finite length \cite[Corollary~6.3.8]{C95}.
Over a finite field $\pi$ is finite-dimensional. Thus all terms
of the restriction filtrations and their subquotients belong to
$\mathcal R_G(B_0)$ or $\mathcal R_G(U)$, respectively. For
\eqref{Bres} we group the irreducible constituents of $[\sigma]$
by their central characters and use exactness of parabolic induction.
On each such summand $Z_A$ acts by the restriction of the Levi
central character.

The Bruhat filtration of $\pi|_{B_0}$ has, on the cell represented by
$w$, the quotient
\[
   \cind_{S_w}^{B_0} W_w,\qquad S_w=P^w\cap B_0,\qquad
   W_w=\sigma^{w^{-1}}|_{D_w},\quad D_w=A(L^w\cap U).
\]
Indeed, the quotient consists of sections compactly supported on
$S_w\backslash B_0$. The coefficient is inflated along the Levi
projection $p_w:S_w\to D_w$ with kernel $V^w\cap U$.

Fix a Bruhat cell indexed by $w$. Put $J_w=S_w\cap U$ and
$U_{L,w}=L^w\cap U$. For a smooth $D_w$-module $E$ write
$T_w(E)=\cind_{S_w}^{B_0}(E\circ p_w)$. By the Bruhat-cell
observation of Section~\ref{sec:grothendieck-B}, $T_w$ takes
$\mathcal R_{L^w}(D_w)$ exactly into $\mathcal R_G(B_0)$:
if $E$ is a subquotient of $\rho|_{D_w}$, its image is a
subquotient of the same Bruhat cell for the parabolic induction
of $\rho$ (transported to $L$). As $p_w$ is the identity on
$Z_A$, this preserves the required central action.

For Conjecture~\ref{conj:restriction-A}, apply $T_w$ to the
asserted relation for $W_w$, with the Levi central character fixed.
A term $\cind_J^{D_w}\chi$, with $A\subseteq J\subseteq D_w$,
therefore contributes, by compact induction in stages,
\[
                  \cind_{p_w^{-1}(J)}^{B_0}(\chi\circ p_w).
\]
The inverse image contains $A$, and the character on $Z_A$ is preserved.

For Conjecture~\ref{conj:restriction-B}, use the Levi relation on
$U_{L,w}$ directly. Inflate a smooth $U_{L,w}$-module $E$ to $J_w$
along $p_w$ and then compactly induce it to $U$. Their composite is exact and takes
$\mathcal R_{L^w}(U_{L,w})$ to $\mathcal R_G(U)$ by the same
Bruhat-cell argument for restriction to $U$.
Thus a term $\cind_J^{L^w\cap U}\one$ contributes
$\cind_{(p_w|_{S_w\cap U})^{-1}(J)}^U\one$ by compact induction in
stages. This inverse image is normalized by $A$ whenever $J$ is.
All filtration terms lie in the appropriate category as submodules
of the ambient restriction. Summing over the Bruhat filtration proves the
two assertions. If the representations appearing in a virtual representation of the Levi subgroup satisfy the corresponding restriction identity, additivity gives the identity for its parabolic induction. More generally, inflation and compact induction are exact and hence induce homomorphisms on Grothendieck groups. Applying these homomorphisms to a restriction identity for the virtual representation itself and summing over the Bruhat cells gives the induced identity.
\end{proof}

\Needspace{5\baselineskip}
\begin{corollary}\label{cor:parabolic-restriction}
Let $\mathbf P=\mathbf L\mathbf V$ be an $F$-parabolic subgroup of $\mathbf G$, and let $\sigma$ be a finite-length smooth admissible representation of $L=\mathbf L(F)$. Choose a minimal $F$-parabolic of $\mathbf L$ and write $U_L$ for the $F$-points of its unipotent radical. Then $\Ind_P^G(\sigma)$ satisfies the restriction identities \eqref{Bres} and \eqref{eq:conjecture-B} in either of the following cases.
\begin{enumerate}
\item $\sigma$ is one-dimensional and trivial on $U_L$.
\item As an $F$-group, $\mathbf L$ is a direct product of an $F$-torus and factors of the form $\GL_{n_i}$ or $\mathrm{PGL}_{n_i}$.
\end{enumerate}
\end{corollary}

\begin{proof}
By additivity on irreducible constituents, taken separately for each
central character, it suffices to treat irreducible $\sigma$.
Proposition~\ref{prop:induction-permanence} then reduces the assertion
to the two identities on $L$. Write $B_{0,L}=AU_L$, with $A$ the group of $F$-points of a maximal split torus in the chosen minimal parabolic. In case (1), triviality on $U_L$ gives $\sigma(au)=\chi(a)$ for
$a\in A$ and $u\in U_L$, where $\chi=\sigma|_A$. Thus
$\sigma|_{B_{0,L}}$ is the inflation of $\chi$ and $\sigma|_{U_L}=\one$.
These are the single terms in \eqref{Bres} and \eqref{eq:conjecture-B}
obtained by taking $H=B_{0,L}$ and $V=U_L$, respectively.
Over a local field the hypothesis on $U_L$ is automatic for a smooth
character of $L$, since a suitable element of $A$ contracts $U_L$ and
$\ker\sigma$ is open. Over a finite field a character of $L$ can be nontrivial
on $U_L$ only in a few exceptional cases arising from small almost simple
factors \cite[Theorem~24.17]{MT11}.

In case (2), apply Theorem~\ref{thm:GL} to each general linear factor. For a projective general linear factor, inflate to $\GL_{n_i}(F)$. The proof of Theorem~\ref{thm:GL} gives its identity in the category with fixed central character. Here this character is trivial. The exact sequences in the proof of Theorem~\ref{thm:GL} have ambient representations with this central character, so those representations descend to $\mathrm{PGL}_{n_i}(F)$, as do the restriction relations. Taking external tensor products, including the character of the torus factor, gives the result for irreducible $\sigma$, and additivity gives the finite-length case.
\end{proof}

\subsection{What remains for Deligne--Lusztig induction}
Let $F$ be finite,  $\mathbf L$ be an $F$-rational Levi subgroup in $\mathbf G$, and $\sigma$ be a representation of $\mathbf L(F)$. Note that $\mathbf L$ may not be a Levi subgroup of an $F$-rational parabolic subgroup, and in this case the Mackey formula for Deligne-Lusztig induction concerns the
composition 
\[
{}^*R_{\mathbf T}^{\mathbf G}  \circ R_{\mathbf L}^{\mathbf G}(\sigma) = \left(R_{\mathbf L}^{\mathbf G}(\sigma) \right)^U
\]
as a representation of $T$, where ${}^*R_{\mathbf T}^{\mathbf G}$ is the {\it Lusztig restriction}.
A more general Mackey formula (often called the \emph{geometric lemma}) for the
restriction of Deligne--Lusztig or Lusztig induction to a parabolic subgroup
seems to be unavailable. Such a formula, even in a less precise form than
its counterpart for ordinary parabolic induction, is what is called for here.  

For Deligne--Lusztig induction, we can reduce the restriction conjectures to the case of elliptic Green functions as follows.


\begin{corollary}\label{cor:elliptic-green}
Let $F$ be a finite field. The following are equivalent, with
$\mathbf G$ ranging over \emph{all} connected reductive $F$-groups:
\begin{enumerate}
\item Both restriction conjectures hold for every uniform representation.
\item For every elliptic maximal $F$-torus $\mathbf T'$ of $\mathbf G$,
the Green function $Q_{\mathbf T'}^{\mathbf G}|_U$ is a linear
combination of $\Ind_V^U\one$ with $V\subset U$ normalized by $A$.
\end{enumerate}
\end{corollary}

\begin{proof}
On unipotent elements, $R_{\mathbf T'}^{\mathbf G}(\theta)$ equals
$Q_{\mathbf T'}^{\mathbf G}$ independently of $\theta$
\cite[Theorem~4.2]{DL76}; hence (1) implies (2).

Conversely, assume (2). Induction on the semisimple rank proves
Conjecture~\ref{conj:restriction-B} for every uniform representation:
a nonelliptic rational maximal torus is contained in a Levi subgroup
of a proper rational parabolic, and one applies transitivity of
Lusztig induction \cite[Proposition~9.1.8]{DM20} and
Proposition~\ref{prop:induction-permanence}.

It remains to prove Conjecture~\ref{conj:restriction-A}. Write
$f=R_{\mathbf T'}^{\mathbf G}(\theta)$ and fix $t\in A$. Put
$\mathbf G_t=C_{\mathbf G}(t)^\circ$ and $U_t=C_U(t)$.
Then $A$ is a maximal split torus of $\mathbf G_t$, and $U_t$ is
the group of rational points of the unipotent radical of its
minimal parabolic subgroup determined by $B$.
The character formula \cite[Theorem~4.2]{DL76} expresses
$u\mapsto f(tu)$ on $U_t$ as a linear combination of Green
functions of $\mathbf G_t$. The already established assertion for
$U_t$ therefore gives
\[
 f(tu)=\sum_V a_{t,V}\,\Ind_V^{U_t}\one(u),\qquad u\in U_t,
\]
where each $V\subseteq U_t$ is normalized by $A$.

For such $V$, identify $(AV)\backslash B_0$ with $V\backslash U$.
The orders of $t$ and $u\in U_t$ are coprime, so a coset is fixed
by $tu$ if and only if it is fixed by both $t$ and $u$.
A coset $Vx$ is fixed by $t$ precisely when
$xtx^{-1}\in Vt$. Since $t$ centralizes $V$, the only element
of $Vt$ of order prime to $p$ is $t$. Thus $x\in U_t$, and
the $t$-fixed cosets are exactly $V\backslash U_t$. Consequently,
\begin{equation}\label{eq:centralizer-permutation}
 \Ind_{AV}^{B_0}\one(tu)=\Ind_V^{U_t}\one(u),\qquad u\in U_t.
\end{equation}

Every element of $tU$ is $U$-conjugate to an element of $tU_t$.
Indeed, if $g\in tU$ has commuting semisimple and unipotent
parts $g_s$ and $g_u$, then $\langle g_s\rangle$ and
$\langle t\rangle$ are complements to $U$ in
$U\rtimes\langle t\rangle$. Conjugacy of complements gives
$x\in U$ with $xg_sx^{-1}=t$, and then
$xg_ux^{-1}\in U_t$. Hence \eqref{eq:centralizer-permutation}
and the expansion of $f(tu)$ determine $f$ on the whole coset
$tU$.

Inflate characters of $A$ to $B_0$ and let $e_t$ be the
characteristic function of $tU$. Fourier inversion on $A$ gives
\[
 e_t(b)=\frac1{|A|}\sum_{\chi\in\widehat A}
              \chi(t)^{-1}\chi(b),\qquad b\in B_0,
\]
where $\widehat A$ is the Pontryagin dual of $A$. 
It follows that
\[
 f|_{B_0}=\frac1{|A|}\sum_{t\in A}\sum_V
       \sum_{\chi\in\widehat A}
       a_{t,V}\chi(t)^{-1}\Ind_{AV}^{B_0}\chi.
\]
Projecting to a prescribed $Z_A$-character retains exactly the
terms with the corresponding restriction $\chi|_{Z_A}$.
Linearity proves Conjecture~\ref{conj:restriction-A} for every
uniform representation.
\end{proof}

We now give further evidence for Conjectures~\ref{conj:restriction-A} and~\ref{conj:restriction-B}.

\subsection{The choice of subgroups}
The Steinberg and principal-series filtrations of
Sections~\ref{sec:steinberg} and~\ref{sec:principal} may suggest that, in the split case, the subgroups
\[
 B_w=B\cap w^{-1}Bw,\qquad w\in W,
\]
are adequate as inducing subgroups in \eqref{Bres}. However, this is not always the case.
The following table summarizes the results of the preceding sections and
calculations done with ChatGPT for the low rank finite groups
${\rm Sp}_4(\mathbb{F}_q)$, ${\rm Sp}_6(\mathbb{F}_q)$ and $G_2(\mathbb{F}_q)$.
Cases outside the stated ranges, such as even $q$, were not calculated;
no different outcome is suggested for those cases. 

\par\addvspace{\baselineskip}
\begin{center}
\renewcommand{\arraystretch}{1.15}
\begin{tabular}{p{.20\textwidth}p{.30\textwidth}p{.40\textwidth}}
\hline
Groups & Range and representations & Adequate inducing subgroups $H_i$ and coefficients $c_i$\\
\hline
$\GL_n(F)$ & Finite or local $F$; all irreducible smooth representations
&  $B_w$; integral coefficients   \\
$\mathrm{Sp}_4(\mathbb{F}_q),\ \mathrm{Sp}_6(\mathbb{F}_q)$ & Every odd $q$; all Deligne--Lusztig characters
&  $B_w$; integral coefficients   \\
$G_2(\mathbb{F}_q)$ & $q=p^f$, $p>3$; all Deligne--Lusztig characters
& Twelve $B_w$ with one additional subgroup $H$; integral coefficients  \\
\hline
\end{tabular}
\end{center}
\par\addvspace{\baselineskip}

For $G_2$, let $a,b$ be respectively the short and long simple roots,
and put $c=a+b$, $d=2a+b$, $e=3a+b$, $f=3a+2b$. If $U_R$ denotes the
subgroup generated by the root groups indexed by $R$, the additional
subgroup is
\[
 H =TU_{\{a,e,f\}}.
\]
Its unipotent part is abelian of order $q^3$ and is not of the form $U\cap w^{-1}Uw$. The thirteen patterns are
\[
 \varnothing,\ a,\ b,\ bc,\ ae,\ ade,\ bcf,\ bcdf,\ aef,\
 adef,\ acdef,\ bcdef,\ abcdef.
\]
They suffice for all $p>3$. 
The twelve Weyl intersections alone fail even over $\C$, already for Coxeter and longest-element
Deligne--Lusztig characters. Thus the general conjecture must allow
more subgroups $A\subseteq H\subseteq B_0\subseteq B$ than $B\cap w^{-1}Bw$. 

\subsection{The rank 3 finite unitary group}
Let $G=\mathrm U_3(\mathbb{F}_q)$ be the finite unitary group, $U$ a fixed
maximal unipotent subgroup, $Z_U$ the center of $U$, and $\rho$ the cuspidal
unipotent character of $G$ of degree $q(q-1)$. Then $|U|=q^3$, $|Z_U|=q$,
with the following information on the character of $\rho$:
\[
 \rho(1)=q(q-1),\quad \rho(z)=-q\ (z\in Z_U\setminus\{1\}),
 \quad \rho(u)=0\ (u\notin Z_U);
\]
see \cite[Section~2.1, Table~3.1]{O13}. Put
\[
 A=\{t_a:=\operatorname{diag}(a,1,a^{-1}):a\in\mathbb F_q^\times\},
 \qquad B_0=AU.
\]
For odd $q$, use Heisenberg coordinates
$U=\mathbb F_{q^2}\times\mathbb F_{q^2}^{-}$, where $\bar x=x^q$,
$\mathbb F_{q^2}^{-}$ is the trace-zero line, and
\[
 (x,s)(y,t)=(x+y,s+t+(y\bar x-x\bar y)/2).
\]
Every $\mathbb F_q$-line $L\subset\mathbb F_{q^2}$ gives a
subgroup $V_L=L\times\{0\}$ which is normalized by $A$. Calculating the characters of the induced representations by counting conjugates, we get:
\begin{equation}\label{eq:unitary-odd-U}
 \rho|_U=\Ind_{V_L}^U\one-\Ind_{Z_UV_L}^U\one.
\end{equation}
Indeed, the difference of the characters of the induced representations above has value $q^2-q$ at $1$, value $-q$ on
$Z_U\setminus\{1\}$, and value $0$ elsewhere, since for a noncentral
$(x,s)$ both induced characters have value $q$ if $x\in L$, and
both vanish otherwise. Similarly, for every $q$ we find that
\begin{equation}\label{eq:unitary-rational}
 \rho|_U=\frac1q\left(\Ind_{\{1\}}^U\one-\Ind_{Z_U}^U\one\right).
\end{equation}
The two subgroups of $U$ used for induction in \eqref{eq:unitary-rational} are normalized by the full rational torus $T=Z_G(A)$. In contrast, the rational maximal torus
$T\subset B$ generally does not normalize $V_L$. This illustrates the
usefulness of $A$ for the integral identity \eqref{eq:unitary-odd-U}. For us,
the essential point is that the subgroups of $U$ used in this integral
identity are normalized by $A$, although they are generally not
normalized by the full maximal torus normalizing $U$.

There are corresponding identities for restriction to $B_0$.
Suppose first that $q$ is odd. Choose a character
$\chi:A\to\C^\times$ with $\chi(t_{-1})=-1$; such a character exists
because $A\simeq\mathbb F_q^\times$ is cyclic of even order.
Extend $\chi$ to $AV_L$ and $AZ_UV_L$ by declaring it to be trivial on $V_L$ and $Z_UV_L$, respectively. Then
\begin{equation}\label{eq:unitary-odd-A}
 \rho|_{B_0}=\Ind_{AV_L}^{B_0}\chi-\Ind_{AZ_UV_L}^{B_0}\chi.
\end{equation}
We omit the proof.

For even $q$, we have instead
\begin{equation}\label{eq:unitary-even-A}
 \rho|_{B_0}=\frac1q\left(\Ind_A^{B_0}\one-\Ind_{AZ_U}^{B_0}\one\right).
\end{equation}

\section*{Declaration on AI usage}
The authors have extensively used ChatGPT in all possible ways: for checking with examples and counter-examples, formulation of results, even completing proofs, besides writing of the paper and locating references. The authors have independently checked all material obtained with the assistance of ChatGPT and take full responsibility for the correctness of the paper.


\begin{thebibliography}{99}

\bibitem[ABS22]{ABS22}
A. Aizenbud, J. Bernstein, and E. Sayag, \emph{Strong density of spherical characters attached to unipotent subgroups}, \href{https://arxiv.org/abs/2202.04984}{arXiv:2202.04984} (2022).

\bibitem[B91]{B91}
A. Borel, \emph{Linear Algebraic Groups}, second edition, Graduate Texts in Mathematics \textbf{126}, Springer-Verlag, New York, 1991.

\bibitem[BDT20]{BDT20}
O. Brunat, O. Dudas, and J. Taylor, \emph{Unitriangular shape of decomposition matrices of unipotent blocks}, Annals of Mathematics (2) \textbf{192} (2020), no.~2, 583--663.

\bibitem[BMRS24]{BMRS24}
M. Bate, B. Martin, G. R\"ohrle, and D. Sercombe, \emph{The subgroup structure of pseudo-reductive groups}, Oberwolfach Preprints 2024-08 (2024), 
 \href{https://doi.org/10.14760/OWP-2024-08}{doi:10.14760/OWP-2024-08}.

\bibitem[BR06]{BR06}
C. Bonnaf\'e and R. Rouquier, \emph{Coxeter orbits and modular
representations}, Nagoya Math. J. \textbf{183} (2006), 1--34.

\bibitem[BS76]{BS76}
A. Borel and J.-P. Serre, \emph{Cohomologie d'immeubles et de groupes $S$-arithm\'etiques},
Topology \textbf{15} (1976), no.~3, 211--232.

\bibitem[BT65]{BT65}
A. Borel and J. Tits, \emph{Groupes r\'eductifs}, Publications Math\'ematiques de l'IH\'ES \textbf{27} (1965), 55--151.

\bibitem[BZ76]{BZ76}
I. N. Bernstein and A. V. Zelevinsky, \emph{Representations of the group $\GL(n,F)$, where $F$ is a local non-Archimedean field}, Russian Mathematical Surveys \textbf{31} (1976), no.~3, 1--68.

\bibitem[BZ77]{BZ77}
I. N. Bernstein and A. V. Zelevinsky, \emph{Induced representations of reductive $p$-adic groups. I}, Annales scientifiques de l'\'Ecole Normale Sup\'erieure (4) \textbf{10} (1977), no.~4, 441--472.

\bibitem[Bo02]{Bo02}
N. Bourbaki, \emph{Lie Groups and Lie Algebras, Chapters 4--6}, Springer, 2002.

\bibitem[C95]{C95}
W. Casselman, \emph{Introduction to the Theory of Admissible Representations of $p$-adic Reductive Groups}, draft, 1 May 1995.

\bibitem[ChevieJL]{ChevieJL}
J. Michel, \emph{Chevie.jl documentation: Unipotent classes of reductive groups},
\href{https://jmichel7.github.io/Chevie.jl/stable/ucl/#Chevie.Ucl.GreenTable}{GreenTable} (split $G_2$ example).

\bibitem[DL76]{DL76}
P. Deligne and G. Lusztig, \emph{Representations of reductive groups over finite fields}, Annals of Mathematics (2) \textbf{103} (1976), 103--161.

\bibitem[DM85]{DM85}
F. Digne and J. Michel, \emph{Fonctions $L$ des vari\'et\'es de Deligne--Lusztig et descente de Shintani}, M\'emoires de la Soci\'et\'e Math\'ematique de France (N.S.) \textbf{20} (1985).

\bibitem[DM20]{DM20}
F. Digne and J. Michel, \emph{Representations of Finite Groups of Lie Type}, second edition, London Mathematical Society Student Texts \textbf{96}, Cambridge University Press, 2020.

\bibitem[F95]{F95}
G. B. Folland, \emph{A Course in Abstract Harmonic Analysis},
Studies in Advanced Mathematics, CRC Press, Boca Raton, FL, 1995.

\bibitem[Gel70]{Gel70}
S. I. Gel'fand, \emph{Representations of the full linear group over a finite field}, Math. USSR-Sb. \textbf{12} (1970), 13--39.

\bibitem[GGM14]{GGM14}
O. Gabber, P. Gille, and L. Moret-Bailly, \emph{Fibr\'es principaux sur les corps valu\'es hens\'eliens}, Algebraic Geometry \textbf{1} (2014), no.~5, 573--612.

\bibitem[GGS17]{GGS17}
R. Gomez, D. Gourevitch, and S. Sahi, \emph{Generalized and degenerate Whittaker models}, Compositio Math. \textbf{153} (2017), no.~2, 223--256.

\bibitem[GRS92]{GRS92}
L. Gurvits, L. Rodman, and T. Shalom, \emph{Controllability and completion of partial upper triangular matrices over rings}, Linear Algebra Appl. \textbf{172} (1992), 135--149.

\bibitem[Gi02]{Gi02}
P. Gille, \emph{Unipotent subgroups of reductive groups in characteristic $p>0$}, Duke Mathematical Journal \textbf{114} (2002), no.~2, 307--328.

\bibitem[Gr02]{Gr02}
J. Grodal, \emph{Higher limits via subgroup complexes}, Annals of Mathematics (2) \textbf{155} (2002), no.~2, 405--457.

\bibitem[K85]{K85}
N. Kawanaka, \emph{Generalized Gelfand--Graev representations and Ennola duality}, in Algebraic Groups and Related Topics (Kyoto/Nagoya, 1983), Advanced Studies in Pure Mathematics \textbf{6}, North-Holland, Amsterdam, 1985, 175--206.

\bibitem[Laz54]{Laz54}
M. Lazard, \emph{Sur les groupes nilpotents et les anneaux de Lie},
Ann. Sci. \'Ecole Norm. Sup. (3) \textbf{71} (1954), no.~2, 101--190.

\bibitem[L92]{L92}
G. Lusztig, \emph{A unipotent support for irreducible representations}, Advances in Mathematics \textbf{94} (1992), no.~2, 139--179.

\bibitem[MT11]{MT11}
G. Malle and D. Testerman, \emph{Linear Algebraic Groups and Finite Groups
of Lie Type}, Cambridge Studies in Advanced Mathematics \textbf{133},
Cambridge University Press, 2011.

\bibitem[MW87]{MW87}
C. M\oe glin and J.-L. Waldspurger, \emph{Mod\`eles de Whittaker d\'eg\'en\'er\'es pour des groupes $p$-adiques}, Mathematische Zeitschrift \textbf{196} (1987), no.~3, 427--452.

\bibitem[NP]{NP}
D. Prasad and N. Sanat, \emph{On the restriction of cuspidal representations to unipotent elements},
Math. Proc. Cambridge Philos. Soc. \textbf{132} (2002), no.~1, 35--56,
\href{https://doi.org/10.1017/S0305004101005527}{doi:10.1017/S0305004101005527}.

\bibitem[O13]{O13}
S. Yu. Orevkov, \emph{Products of conjugacy classes in finite unitary groups $GU(3,q^2)$ and $SU(3,q^2)$}, Annales de la Facult\'e des Sciences de Toulouse (6) \textbf{22} (2013), 219--251.

\bibitem[R75]{R75}
F. Rodier, \emph{Mod\`eles de Whittaker des repr\'esentations admissibles
des groupes r\'eductifs $p$-adiques quasi-d\'eploy\'es}, manuscript,
revised 1975, available at
\url{https://www.normalesup.org/~rodier/Whittaker.pdf}.

\bibitem[T16]{T16}
J. Taylor, \emph{Generalized Gelfand--Graev representations in small characteristics}, Nagoya Mathematical Journal \textbf{224} (2016), 93--167.

\end{thebibliography}
\end{document}